\documentclass[12pt,a4paper]{elsarticle}

\usepackage[dvipsnames]{xcolor}
\usepackage{graphicx}
\usepackage{amssymb}
\usepackage{amsmath}
\allowdisplaybreaks
\usepackage{bm}
\usepackage{subfigure}
\usepackage[T1]{fontenc}
\usepackage{amsthm}
\usepackage{tikz}
\usepackage{color}
\usepackage{MnSymbol}
\usepackage{bbding}

\usetikzlibrary{external}
\usetikzlibrary{shapes}

\usepackage{hyperref}
\usepackage{pgfplots}

\usepackage{etoolbox}
\usetikzlibrary{pgfplots.patchplots}
\usetikzlibrary{pgfplots.colormaps}
\usetikzlibrary{pgfplots.groupplots}

\usepackage{algorithm2e}
\usepackage{soul}
\usepackage{enumerate}
\usepackage{multirow}
\usepackage{colortbl}
\usepackage{hhline}

\usepackage[title,toc,page]{appendix}

\usetikzlibrary{positioning,arrows,mindmap,backgrounds, shapes, decorations.fractals, calc, intersections,decorations.pathmorphing,patterns}

\usepackage[colorinlistoftodos, obeyFinal]{todonotes}

\makeatletter
\newcommand{\vast}{\bBigg@{17.5}}
\newcommand{\Vast}{\bBigg@{9}}
\newcommand{\thickhline}{
	\noalign {\ifnum 0=`}\fi \hrule height 1pt
	\futurelet \reserved@a \@xhline
}
\newcolumntype{?}{!{\vrule width 1pt}}
\makeatother
\newcommand\Tstrut{\rule{0pt}{2.6ex}}         
\newcommand\Bstrut{\rule[-0.9ex]{0pt}{0pt}}   

\newcommand{\markerone}{\raisebox{0.0 cm}{\tikz{\node[draw=blue,scale=0.7,circle,thick](){};}}}
\newcommand{\markertwo}{\raisebox{0.0cm}{\tikz{\node[draw=green,scale=0.5,diamond,fill=green,thick](){};}}}
\newcommand{\markerthree}{\raisebox{0.0cm}{\tikz{\node[draw=red,scale=0.5,cross out,rotate=45,fill=red,thick](){};}}}
\newcommand{\markerfour}{\raisebox{0.5pt}{\tikz{\node[draw,scale=0.4,cross out,fill=black,thick](){};}}}
\newcommand{\markerfive}{\raisebox{0.0cm}{\tikz{\node[draw=red,line width=0.8pt,scale=0.4,circle,fill=none](){};}}}

\newtheorem{remark}{Remark}

\newcommand{\eg}{{e.g.,~}}

\def\rmd{{\mathrm{d}}}

\newcommand*{\Scale}[2][4]{\scalebox{#1}{$#2$}}%

\def\be{\begin{equation}}
\def\ee{\end{equation}}
\def\ba{\begin{array}}
\def\ea{\end{array}}

\definecolor{DarkScarlet}{rgb}{0.34, 0.01, 0.1}
\definecolor{dgreen}{RGB}{0,139,0}

\title{Fast and accurate elastic analysis of laminated composite plates via isogeometric collocation and an equilibrium-based stress recovery approach}
		\author[pavia1]{Alessia Patton\corref{cor2}}
	\author[UTA]{John-Eric Dufour}
	\author[lausanne]{Pablo Antolin}
	\author[pavia1,imati]{Alessandro Reali}	
	\cortext[cor2]{Corresponding author. E-mail: \textsf{alessia.patton01@universitadipavia.it}}
	\address[pavia1]{Department of Civil Engineering and Architecture - University of Pavia\\ Via Ferrata, 3, 27100 Pavia, Italy}
	
	\address[UTA]{Mechanical and Aerospace Engineering Department\\
		University of Texas at Arlington\\
		500 W 1st St, Arlington, TX 76010}
	
	\address[lausanne]{Institute of Mathematics - \'Ecole Polytechnique F\'ed\'erale de Lausanne\\ CH-1015 Lausanne, Switzerland}
	\date{\empty}

\address[imati]{Instituto di Matematica Applicata e Tecnologie Informatiche ``E.~Magenes'' (CNR)\\
via Ferrata 1, 27100 Pavia, Italy}

\pgfplotsset{stress plot style/.style={
  width=.48\textwidth,
  y label style={at={(axis description cs:0.15,.5)},anchor=south},
  ylabel={$x_3$ [mm]}}}

\newif\ifrecompiletikz
 \recompiletikztrue

\ifrecompiletikz
 \usetikzlibrary{external}
 \tikzexternaldisable
\fi

\begin{document}
		
\begin{abstract}
A novel approach which combines isogeometric collocation and an equilibrium-based stress recovery technique is applied to analyze laminated composite plates. Isogeometric collocation is an appealing strong form alternative to standard Galerkin approaches, able to achieve high order convergence rates coupled with a significantly reduced computational cost.
Laminated composite plates are herein conveniently modeled considering only one element through the thickness with homogenized material properties. This guarantees accurate results in terms of displacements and in-plane stress components.
To recover an accurate out-of-plane stress state, equilibrium is imposed in strong form as a post-processing correction step, which requires the shape functions to be highly continuous.
This continuity demand is fully granted by isogeometric analysis properties, and excellent results are obtained using a minimal number of collocation points per direction, particularly for increasing values of length-to-thickness plate ratio and number of layers.
\end{abstract}


\setcounter{tocdepth}{1}
\maketitle
\section{Introduction}\label{sec:intro}
\begin{figure}[!htbp]
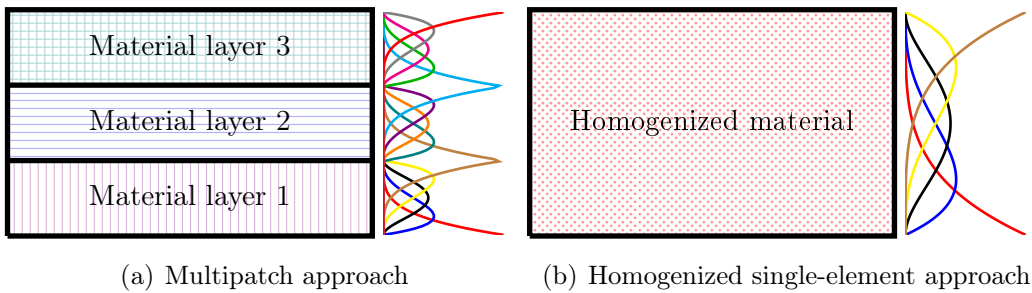

	\centering
	\subfigure[Multipatch approach\label{subfig-1:layerwise}]{\ifrecompiletikz\tikzsetnextfilename{fig_01_a}\tikzexternalenable\begin{tikzpicture}
\def\L{4.8};
\def\t{3};
\def\nlayers{3};
\def\dt{\t/\nlayers};
\def\nq{4};
\def\ddt{\dt/\nq};
\def\tol{0.05};
\def\arrL{1.1}

horizontal lines,
vertical lines,
north east lines,
north west lines,
grid,
crosshatch,
dots,
crosshatch dots,

\draw[pattern=vertical lines, pattern color=violet!30,draw=none] (0,0*\dt) rectangle (\L,1*\dt);
\draw[pattern=horizontal lines, pattern color=blue!30,draw=none]  (0,1*\dt) rectangle (\L,2*\dt);
\draw[pattern=grid, pattern color=teal!30,draw=none]             (0,2*\dt) rectangle (\L,3*\dt);
\draw[line width=2pt] (0, 0)--(\L, 0)--(\L, \t)--(0, \t)--(0,0);

\foreach \x in {1,...,\nlayers}
{
\def\xx{\dt*\x-\dt/2};
\node at (\L*2/4,\xx) {\small Material layer \x};
}

\foreach \x in {2,...,\nlayers}
{
\def\ti{\dt*\x-\dt};
\draw[line width=2pt](0, \ti) -- (\L, \ti);
}

\node [anchor=south west]
(plot) at (\L, -0.13)
{
\begin{tikzpicture}[trim axis left]
\def\width{90};
\def\height{129};

\begin{axis}
[
   line width=1pt,
   cycle list name=color list,
   width=\width, height=\height,
   xtick={\empty},
   ymin = 0.0,  ymax=1,
   xmin = 0.0,  xmax=1,
   ytick={\empty},
   axis x line=none,
   axis y line=none,
]
\foreach \x in {1,...,13}
{
\addplot table [x index=\x, y index = 0] {Images/functions/c0.dat};}
\end{axis}
\end{tikzpicture}
};
\end{tikzpicture}\tikzexternaldisable\else\includegraphics{fig_01_a.pdf}\fi}\hfill
	\subfigure[Homogenized single-element approach\label{subfig-2:HSE}]{\ifrecompiletikz\tikzsetnextfilename{fig_01_b}\tikzexternalenable\begin{tikzpicture}
\def\L{4.8};
\def\t{3};
\def\nlayers{3};
\def\dt{\t/\nlayers};
\def\nq{4};
\def\ddt{\dt/\nq}

horizontal lines,
vertical lines,
north east lines,
north west lines,
grid,
crosshatch,
dots,
crosshatch dots

\draw[pattern=crosshatch dots, pattern color=red!40,draw=none]   (0,0) rectangle (\L,\t);
\draw[line width=2pt] (0, 0)--(\L, 0)--(\L, \t)--(0, \t)--(0,0);

\node [black] at (\L*0.5,\t*0.5) {\small Homogenized material};
\node [anchor=south west]
(plot) at (\L, -0.13)
{
\begin{tikzpicture}[trim axis left]
\def\width{90};
\def\height{129};

\begin{axis}
[
   line width=1pt,
   cycle list name=color list,
   width=\width, height=\height,
   xtick={\empty},
   ymin = 0.0,  ymax=1,
   xmin = 0.0,  xmax=1,
   ytick={\empty},
   axis x line=none,
   axis y line=none,
]
\foreach \x in {1,...,5}
{
\addplot table [x index=\x, y index = 0] {Images/functions/single_element.dat};
}
\end{axis}
\end{tikzpicture}
};
\end{tikzpicture}\tikzexternaldisable\else\includegraphics{fig_01_b.pdf}\fi}
	\caption{Layerwise approach and homogenized single-element example of isogeometric shape functions for a degree of approximation equal to 4.}
	\label{fig:layerwiseVsHSE}
\end{figure}

Composite materials consist of two or more materials which combined present enhanced properties that could not be acquired employing any of the constituents alone (see, \eg\cite{Gibson1994,Jones1999,Reddy2003,Vinson1986} and references therein). The interest for composite structures in the engineering community has constantly grown in recent years due to their appealing mechanical properties such as increased stiffness and strength, reduced weight, improved corrosion and wear resistance, just to recall some of them. The majority of man-made composite materials consists of reinforced fibers embedded in a base material, called matrix (see, \eg\cite{Hashin1972,Reddy2003}).
The matrix material keeps the fibers together, acts as a load-transfer medium between fibers, process which takes place through shear stresses, and protects those elements from being exposed to the environment, while the resistence properties of composites are given by the fibers which are stiffer and stronger than the soft matrix.
In this paper we focus on laminated composite materials, which are formed by a collection of building blocks or plies, stacked to achieve the desired stiffness and thickness. 
For this kind of structures also simple loading conditions, such as traction or bending, cause a complex 3D stress state because of the difference in the material properties between the layers, which may lead to delamination, consequently requiring an accurate stress evaluation through the thickness (see, \eg\cite{Reddy2003,Sridharan2008}).
As an alternative to two-dimensional theories, often insufficiently accurate to depict delamination and interlaminar damage, and to layerwise theories, which typically show an high computational cost, a novel method combining an isogeometric analysis (IGA) Galerkin approach with a stress-recovery technique has been recently proposed in~\cite{Dufour2018}. 
Introduced in 2005 by Hughes et al.~\cite{Hughes2005}, IGA aims at integrating design and analysis employing shape functions typically belonging to Computer Aided Design field (such as, B-Splines and NURBS). 
Using the same shape functions to approximate both geometry and field variables leads to a cost-saving simplification of expensive mesh generation and refinement processes required by standard finite element analyisis. 
One of the most important features of IGA is the high-regularity of its basis functions leading to superior approximation properties. IGA proved to be successful in a wide variety of problems ranging from solids and structures (see, \eg\cite{Auricchio2010b,Borden2012,Caseiro2014,Dhote2014,Elguedj2014,Hughes2008,Hughes2014,Lipton2010,Morganti2015}) to fluids (see, \eg\cite{Akkerman2008,Buffa2011,Gomez2010,Liu2013}), fluid-structure interaction (see, \eg\cite{Bazilevs2011,Hsu2015}), opening also the door to geometrically flexible discretizations of higher-order partial differential equations in primal form as in~\cite{Auricchio2007,Gomez2008,Kiendl2009,Lorenzo2018}.
However, a well-known important issue of IGA is related to the development of efficient integration rules when higher-order approximations are employed (see, \eg\cite{Auricchio2012b,Fahrendof2018,Hughes2010,Sangalli2018}). In attempt to address this problem taking full advantage of the special possibilities offered by IGA, isogeometric collocation (IGA-C) schemes have been proposed in~\cite{Auricchio2010a}. The aim was to optimize the computational cost still relying on IGA geometrical flexibility and accuracy.
Collocation main idea, in contrast to Galerkin-type formulations, consists in the discretization of the governing partial differential equations in strong form, evaluated at suitable points. Since integration is not required, isogeometric collocation results in a very fast method providing superior performance in terms of accuracy-to-computational effort ratio with respect to Galerkin formulations, in particular when higher-order approximation degrees are adopted (see~\cite{Schillinger2013}). 
Isogeometric collocation has been particularly successful in the context of structural elements, where isogeometric collocation has proven to be particularly stable in the context of mixed methods. In particular, Bernoulli-Euler beam and Kirchhoff plate elements have been proposed~\cite{Reali2015b}, while mixed formulations both for Timoshenko initially-straight planar~\cite{daVeiga2012} and non-prismatic~\cite{Balduzzi2017} beams as well as for curved spatial rods~\cite{Auricchio2013} have been introduced and studied, and then effectively extended to the geometrically nonlinear case~\cite{Kiendl2018,Marino2016,Marino2017,Marino2019,Weeger2017a,Weeger2018}. Isogeometric collocation has been moreover successfully applied to the solution of Reissner-Mindlin plate problems in~\cite{Kiendl2015a}, and new formulations for shear-deformable beams~\cite{Kiendl2015b,Kiendl2018}, as well as shells~\cite{Kiendl2017,Maurin2018} have been solved also via IGA collocation. Since its introduction, many promising significant works on isogeometric collocation methods have been published also in other fields, including phase-field modeling~\cite{Gomez2014}, contact~\cite{DeLorenzis2015,Kruse2015,Weeger2017b} and poromechanics~\cite{Morganti2018}. Moreover,
combinations with different spline spaces, like hierarchical splines, generalized B-Splines, and T-Splines, have been successfully tested in~\cite{Casquero2015,Manni2015,Schillinger2013}, while alternative effective selection strategies for collocation points have been proposed in~\cite{Anitescu2015,Gomez2014,Montardini2017}.  
IGA-Galerkin methods have already been used to solve composite laminate problems, especially relying on high-order theories for enhanced plate and shell theories~\cite{Farzam2018,Kapoor2013,Remmers2015,Shi--Dong2018,Thai2015}. Recently an isogeometric collocation numerical formulation has been proposed~\cite{Pavan2017} to study Reissner-Mindlin composite plates. Other Galerkin methods~\cite{Guo2014,Guo2015,Remmers2015} compute instead a full 3D stress state using isogeometric analysis, applying a layerwise technique.
This method can also be applied to isogeometric collocation adopting a multipatch approach, which models each layer as a patch (see Figure~\ref{subfig-1:layerwise}), enforcing normal stress continuity at the inter-patches boundaries~\cite{Auricchio2012a}. 
Clearly a layerwise method exploits a number of degrees of freedom directly proportional to the number of layers, inevitably leading to high computational costs.
In this paper we apply a single patch 3D isogeometric collocation method to analyze the behavior of composite plates. We adopt a homogenized single-element approach (see Figure~\ref{subfig-2:HSE}), which conveniently uses one element through the thickness, coupled with a post-processing technique in order to recover a proper out-of-plane stress state. This method is significantly less expensive compared to a layerwise approach since employs a considerably lower number of degrees of freedom. 
\newline
The post-processing approach, first proposed in~\cite{Dufour2018}, takes inspiration from recovery techniques which can be found in~\cite{Daghia2008,deMiranda2002,Engblom1985,Kapoor2013,Pryor1971,Ubertini2004} and is based on the direct integration of the equilibrium equations to compute the out-of-plane stress components from the in-plane ones directly derived from a coarse displacement solution.\newline 
The structure of the paper is organized as follows.
In Section~\ref{sec:IGA-Cortho} the fundamental concepts of multivariate B-Splines and NURBS are presented, followed by an introduction to isogeometric collocation and a description of our IGA-C scheme for orthotropic elasticity.
In Section~\ref{sec:IGA-Cstrategies} we define our isogeometric collocation strategy to study laminated plates, which combines a homogenized single-element approach with an equilibrium-based stress recovery technique. 
In Section~\ref{sec:results} we present our reference test case and provide results for the single-element approach.
Several numerical benchmarks are displayed, which show a significant improvement between non-treated and post-processed out-of-plane stress components.
Finally we provide some mesh sensitivity tests considering an increasing length-to-thickness ratio and numbers of layers to show the effectiveness of the method. We draw our conclusions in Section~\ref{sec:conclusions}.
		
\section{Isogeometric Collocation: Basics and application to orthotropic elasticity}\label{sec:IGA-Cortho}
In this section we introduce the notions of multivariate B-Splines and NURBS, provide some details regarding isogeometric collocation and describe our collocation scheme in the context of linear orthotropic elasticity.

\subsection{Multivariate B-Splines and NURBS}\label{subsec:IGAsf}
In the following, we introduce the basic definitions and notations about multivariate B-Splines and NURBS. For further details, readers may refer to~\cite{Cottrell2007,Hughes2005,Piegl1997}, and references therein. Multivariate B-Splines are generated through the tensor product of univariate B-Splines. We denote with $d_{p}$ the dimension of the parametric space and therefore $d_{p}$
univariate knot vectors have to be introduced as
\begin{equation}
\Theta = \{\theta_{1}^{d},...,\theta_{m_{d}+p_{d}+1}^{d}\}\hspace{1cm}d = 1, ..., d_{p}\,,\label{eq:knotvectors}
\end{equation}
where $p_{d}$ represents the polynomial degree in the parametric direction $d$,
and $m_{d}$ is the associated number of basis functions.
Given the univariate basis functions $N^{d}_{i_{d},p_{d}}$ associated to each parametric direction $\xi^{d}$, the multivariate basis functions $B_{\textbf{i},\textbf{p}}(\boldsymbol{\xi})$ are obtained as:
\begin{equation}
B_{\textbf{i},\textbf{p}}(\boldsymbol{\xi})=\prod\limits_{d=1}^{d_{p}}N_{i_{d},p_{d}}(\xi^{d})\,,\label{eq:B-Splines}
\end{equation}
where $\textbf{i} = \{i_{1}, ..., i_{d_{p}}\}$ plays the role of a multi-index which describes the considered position in the tensor product structure, $\textbf{p} = \{{p_{1}, ..., p_{d}}\}$ indicates the polynomial degrees, and
$\boldsymbol{\xi} = \{\xi^{1},...,\xi^{d_{p}}\}$ represents the vector of the parametric coordinates in each parametric direction $d$. B-Spline multidimensional geometries are built from a linear combination of multivariate B-Spline basis functions as follows
\begin{equation}
\textbf{S}(\boldsymbol{\xi})=\sum\limits_{\textbf{i}}B_{\textbf{i},\textbf{p}}(\boldsymbol{\xi})\textbf{P}_\textbf{i}\,,\label{eq:B-Splinemultigeom}
\end{equation}
where the coefficients $\textbf{P}_\textbf{i}\in\mathbb{R}^{d_{s}}$ of the linear combination are the so-called
control points ($d_{s}$ is the dimension of the physical space) and the summation is extended to all combinations of the multi-index \textbf{i}.
NURBS geometries in $\mathbb{R}^{d_{s}}$ are instead obtained from a projective transformation of their B-Spline counterparts in $\mathbb{R}^{d_{s}+1}$.
Defining $w_{\textbf{i}}$ as the collection of weights according to the multi-index \textbf{i}, multivariate NURBS basis functions are obtained as
\begin{equation}
R_{\textbf{i},\textbf{p}}(\boldsymbol{\xi})=\frac{B_{\textbf{i},\textbf{p}}(\boldsymbol{\xi})w_{\textbf{i}}}{\sum_{\textbf{j}}B_{\textbf{j},\textbf{p}}(\boldsymbol{\xi})w_{\textbf{j}}}\label{eq:NURBS}
\end{equation}
and NURBS multidimensional geometries are built as 
\begin{equation}
\textbf{S}(\boldsymbol{\xi})=\sum\limits_{\textbf{i}}R_{\textbf{i},\textbf{p}}(\boldsymbol{\xi})\textbf{P}_\textbf{i}\,.\label{eq:NURBmultigeom}
\end{equation}

\subsection{An Introduction to Isogeometric collocation}\label{subsec:introIGA-C}
Collocation methods have been introduced within isogeometric analysis as an attempt to address a well-known important issue of early IGA-Galerkin formulations, related to the development of efficient integration rules for higher-order approximations. In fact, element-wise Gauss quadrature, typically used for finite elements and originally adopted for Galerkin-based IGA, does
not properly take into account inter-element higher continuity leading to sub-optimal array formation and assembly
costs, significantly affecting the performance of IGA methods. Isogeometric collocation aimed at optimizing computational cost, since it may be viewed as a variant of one-point quadrature numerical scheme, still taking advantage of IGA geometrical flexibility and accuracy.
Collocation methods are based on the direct discretization in strong form of the differential equations governing the problem evaluated at suitable points.
The isoparametric paradigm is adopted and the same basis functions are used to
describe both geometry and problem unknowns. Once the approximations are carried out, as in a typical Galerkin-IGA context, by means of a linear combinations of IGA basis functions and control variables, the discrete differential equations are collocated at each collocation point. 
Consequently a delicate issue is represented by the determination of suitable collocation points. A widespread approach which is proposed in
the engineering literature is to collocate at the images of Greville abscissae (see, \eg\cite{Johnson2005}), but this represents just the simplest possible option (see, \eg\cite{deBoor1973,Demko1985} for alternative choices).
Along each parametric direction $d$, Greville abscissae consist of a set of $m^d$ points, obtained from the knot vector components, $\theta^d_{i}$, as 
\begin{equation}
\overline{\theta}^d_{i}=\frac{\theta^d_{i+1}+\theta^d_{i+2}+...+\theta^d_{i+p}}{p_d}\hspace{1cm}i = 1,...,m^d\,,\label{eq:greville}
\end{equation}
$p_d$ being the degree of approximation.
Since the approximation is performed through direct collocation of the differential equations, no integrals need to be computed and consequently, evaluation and assembly operations lead to a significantly reduced computational cost.

\subsection{Numerical formulation for orthotropic elasticity}\label{subsec:numericalIGA-C}
Once a strategy to select collocation points and compute IGA shape functions is set, a proper description of the  equations in strong form for the problem under examination is required, as mentioned in Section~\ref{subsec:introIGA-C}.
We therefore recall the classical elasticity problem in strong form considering a small strain regime and detail equilibrium equations using Einstein's notation \eqref{eq:equilibrium}. 
The following notations are used:
$\Omega\subset\mathbb{R}^3$, is an open bounded domain, representing
an elastic three-dimensional body, $\Gamma_{N}$ and $\Gamma_{D}$ are defined as boundary portions subjected respectively to Neumann and Dirichlet conditions such that $\Gamma_{N}\cup\Gamma_{D}=\partial\Omega$ and $\Gamma_{N}\cap\Gamma_{D}=\emptyset$. Accordingly, the equilibrium equations and the corresponding boundary conditions are: 
\begin{subequations}
\begin{align}
&\sigma_{ij,j}+b_{i}=0\hspace{2.6cm}\text{in}\;\Omega\label{eq:equilibriumint}\\
&\sigma_{ij}n_{j}=t_{i}\hspace{3.05cm}\text{on}\;\Gamma_{N}\label{eq:equilibriumNeumann}\\
&u_{i}=\overline{u}_{i}\hspace{3.5cm}\text{on}\;\Gamma_{D}\label{eq:equilibriumDirichlet}
\end{align}\label{eq:equilibrium}%
\end{subequations}
where $\sigma_{ij}$ and $u_{i}$ represent respectively the Cauchy stress and displacement components, while $b_{i}$ and $t_{i}$ the volume and traction forces, $n_{j}$ the outward normal, and $\overline{u}_{i}$ the prescribed displacements.
The elasticity problem is finally completed by the kinematic relations in small strain
\begin{equation}
\varepsilon_{ij}=\frac{u_{i,j}+u_{j,i}}{2}\,,\\\label{eq:kinematics}
\end{equation}
as well as by the constitutive equations
\begin{equation}
\sigma_{ij}=\mathbb{C}_{ijkm}{\varepsilon}_{km}\,,\label{eq:constlaw}
\end{equation}
where $\mathbb{C}_{ijkm}$ is the fourth order elasicity tensor.\newline
As we described in Section~\ref{sec:intro}, the basic building block of a laminate is a lamina, i.e., a flat arrangement of unidirectional  fibers,  considering the simplest case, embedded in a matrix.
In order to increase the composite resistance properties cross-ply laminates can be employed (i.e.,  all the plies used to form the composite stacking sequence are piled alternating different fiber layers orientations) in which all unidirectional layers are individually orthotropic.
Since the proposed collocation approach uses one element through the thickness to model the composite plate as a homogenized single building block, we focus in this section on the collocation formulation for a plate formed by only one orthotropic elastic lamina.
Considering three mutually orthogonal planes of material symmetry for each ply, the number of elastic coefficients of the fourth order elasticity tensor $\mathbb{C}_{ijkm}$ is reduced to 9 in Voigt notation, that can be expressed in terms of engineering constants as
\newline\begin{equation}\Scale[0.88]{
\mathbb{C}=\renewcommand\arraystretch{1.75}\begin{bmatrix}
\mathbb{C}_{11} & \mathbb{C}_{12} & \mathbb{C}_{13} & 0 & 0 & 0\\
 & \mathbb{C}_{22} & \mathbb{C}_{23} & 0 & 0 & 0\\
 &  & \mathbb{C}_{33} & 0 & 0 & 0\\
 &  symm&  & \mathbb{C}_{44} & 0 & 0\\
 &  &  &  & \mathbb{C}_{55} & 0\\
 &  &  &  &  & \mathbb{C}_{66}\\
\end{bmatrix}=\vast{[}\begin{matrix}
\cfrac{1}{E_{1}} & -\cfrac{\nu_{12}}{E_{1}} & -\cfrac{\nu_{13}}{E_{1}} & 0 & 0 & 0\\
& \cfrac{1}{E_{2}} & -\cfrac{\nu_{23}}{E_{2}} & 0 & 0 & 0\\
 &  & \cfrac{1}{E_{3}} & 0 & 0 & 0\\
 &  symm&  & \cfrac{1}{G_{23}} & 0 & 0\\
 &  &  &  & \cfrac{1}{G_{13}} & 0\\
 &  &  &  &  & \cfrac{1}{G_{12}}\\
\end{matrix}\vast{]}^{-1}\,.}\label{eq:CIVeltensor}
\end{equation}
The displacement field is then approximate as a linear combination of NURBS multivariate shape functions and control points as follows
\begin{equation}
\begin{aligned}
&\textbf{u}(\boldsymbol{\xi})=R_{\textbf{i},\textbf{p}}(\boldsymbol{\xi})\hat{\textbf{u}}_{\textbf{i}}\,,\\
&\textbf{v}(\boldsymbol{\xi})=R_{\textbf{i},\textbf{p}}(\boldsymbol{\xi})\hat{\textbf{v}}_{\textbf{i}}\,,\\
&\textbf{w}(\boldsymbol{\xi})=R_{\textbf{i},\textbf{p}}(\boldsymbol{\xi})\hat{\textbf{w}}_{\textbf{i}}\,.
\end{aligned}\label{eq:displapprox}
\end{equation}
Having defined $\boldsymbol{\tau}$ as the matrix of collocation points, we insert the approximations~\eqref{eq:displapprox} into kinematics equations~\eqref{eq:kinematics} and we combine the obtained  expressions with the constitutive relations~\eqref{eq:constlaw}. Finally we substitute into equilibrium equations~\eqref{eq:equilibriumint} obtaining
\begin{subequations}\label{eq:IGA-Cinternal}
\begin{align}\tag{\ref{eq:IGA-Cinternal}}
&\begin{bmatrix}
\textbf{K}_{11}(\boldsymbol{\tau})&\textbf{K}_{12}(\boldsymbol{\tau})&\textbf{K}_{13}(\boldsymbol{\tau})\\
&\textbf{K}_{22}(\boldsymbol{\tau})&\textbf{K}_{23}(\boldsymbol{\tau})\\
symm&&\textbf{K}_{33}(\boldsymbol{\tau})
\end{bmatrix}\cdot\begin{pmatrix}\hat{\textbf{u}}_{\textbf{i}}\\
\hat{\textbf{v}}_{\textbf{i}}\\
\hat{\textbf{w}}_{\textbf{i}}\end{pmatrix}=-\textbf{b}(\boldsymbol{\tau}),\hspace{1cm}\forall\boldsymbol{\tau}\in\Omega\,,\\
\intertext{where $\textbf{K}_{ij}(\boldsymbol{\tau})$ cofficients can be expressed as}
&\textbf{K}_{11}(\boldsymbol{\tau})=\mathbb{C}_{11}\cfrac{\partial^2{R_{\textbf{i},\textbf{p}}}(\boldsymbol{\tau})}{\partial{x_1}^2}+\mathbb{C}_{66}\cfrac{\partial^2{R_{\textbf{i},\textbf{p}}}(\boldsymbol{\tau})}{\partial{x_2}^2}+\mathbb{C}_{55}\cfrac{\partial^2{R_{\textbf{i},\textbf{p}}}(\boldsymbol{\tau})}{\partial{x_3}^2}\,,\label{eq:IGA-CinternalK11}\\
&\textbf{K}_{22}(\boldsymbol{\tau})=\mathbb{C}_{66}\cfrac{\partial^2{R_{\textbf{i},\textbf{p}}}(\boldsymbol{\tau})}{\partial{x_1}^2}+\mathbb{C}_{22}\cfrac{\partial^2{R_{\textbf{i},\textbf{p}}}(\boldsymbol{\tau})}{\partial{x_2}^2}+\mathbb{C}_{44}\cfrac{\partial^2{R_{\textbf{i},\textbf{p}}}(\boldsymbol{\tau})}{\partial{x_3}^2}\,,\label{eq:IGA-CinternalK22}\\
&\textbf{K}_{33}(\boldsymbol{\tau})=\mathbb{C}_{55}\cfrac{\partial^2{R_{\textbf{i},\textbf{p}}}(\boldsymbol{\tau})}{\partial{x_1}^2}+\mathbb{C}_{44}\cfrac{\partial^2{R_{\textbf{i},\textbf{p}}}(\boldsymbol{\tau})}{\partial{x_2}^2}+\mathbb{C}_{33}\cfrac{\partial^2{R_{\textbf{i},\textbf{p}}}(\boldsymbol{\tau})}{\partial{x_3}^2}\,,\label{eq:IGA-CinternalK33}\\
&\textbf{K}_{23}(\boldsymbol{\tau})=(\mathbb{C}_{23}+\mathbb{C}_{44})\cfrac{\partial^2{R_{\textbf{i},\textbf{p}}}(\boldsymbol{\tau})}{\partial{x_2}\partial{x_3}}\,,\label{eq:IGA-CinternalK23}\\
&\textbf{K}_{13}(\boldsymbol{\tau})=(\mathbb{C}_{13}+\mathbb{C}_{55})\cfrac{\partial^2{R_{\textbf{i},\textbf{p}}}(\boldsymbol{\tau})}{\partial{x_1}\partial{x_3}}\,,\label{eq:IGA-CinternalK31}\\
&\textbf{K}_{12}(\boldsymbol{\tau})=(\mathbb{C}_{12}+\mathbb{C}_{66})\cfrac{\partial^2{R_{\textbf{i},\textbf{p}}}(\boldsymbol{\tau})}{\partial{x_1}\partial{x_2}}\,,\label{eq:IGA-CinternalK12}
\end{align}
\end{subequations}
and substituting in~\eqref{eq:equilibriumNeumann} we obtain:
\begin{subequations}\label{eq:IGA-CNeumann}
\begin{align}\tag{\ref{eq:IGA-CNeumann}}
&\begin{bmatrix}
\tilde{\textbf{K}}_{11}(\boldsymbol{\tau})&\tilde{\textbf{K}}_{12}(\boldsymbol{\tau})&\tilde{\textbf{K}}_{13}(\boldsymbol{\tau})\\
&\tilde{\textbf{K}}_{22}(\boldsymbol{\tau})&\tilde{\textbf{K}}_{23}(\boldsymbol{\tau})\\
symm&&\tilde{\textbf{K}}_{33}(\boldsymbol{\tau})
\end{bmatrix}\cdot\begin{pmatrix}\hat{\textbf{u}}_{\textbf{i}}\\
\hat{\textbf{v}}_{\textbf{i}}\\
\hat{\textbf{w}}_{\textbf{i}}\end{pmatrix}=\textbf{t}(\boldsymbol{\tau}),\hspace{1cm}\forall\boldsymbol{\tau}\in\Gamma_{N}\\
\intertext{with $\tilde{\textbf{K}}_{ij}(\boldsymbol{\tau})$ components having the following form}
&\tilde{\textbf{K}}_{11}(\boldsymbol{\tau})=\mathbb{C}_{11}\cfrac{\partial{R_{\textbf{i},\textbf{p}}}(\boldsymbol{\tau})}{\partial{x_1}}n_{1}+\mathbb{C}_{66}\cfrac{\partial{R_{\textbf{i},\textbf{p}}}(\boldsymbol{\tau})}{\partial{x_2}}n_{2}+\mathbb{C}_{55}\cfrac{\partial{R_{\textbf{i},\textbf{p}}}(\boldsymbol{\tau})}{\partial{x_3}}n_{3}\,,\label{eq:IGA-CNeumannK11}\\
&\tilde{\textbf{K}}_{22}(\boldsymbol{\tau})=\mathbb{C}_{66}\cfrac{\partial{R_{\textbf{i},\textbf{p}}}(\boldsymbol{\tau})}{\partial{x_1}}n_{1}+\mathbb{C}_{22}\cfrac{\partial{R_{\textbf{i},\textbf{p}}}(\boldsymbol{\tau})}{\partial{x_2}}n_{2}+\mathbb{C}_{44}\cfrac{\partial{R_{\textbf{i},\textbf{p}}}(\boldsymbol{\tau})}{\partial{x_3}}n_{3}\,,\label{eq:IGA-CNeumannK22}\\
&\tilde{\textbf{K}}_{33}(\boldsymbol{\tau})=\mathbb{C}_{55}\cfrac{\partial{R_{\textbf{i},\textbf{p}}}(\boldsymbol{\tau})}{\partial{x_1}}n_{1}+\mathbb{C}_{44}\cfrac{\partial{R_{\textbf{i},\textbf{p}}}(\boldsymbol{\tau})}{\partial{x_2}}n_{2}+\mathbb{C}_{33}\cfrac{\partial{R_{\textbf{i},\textbf{p}}}(\boldsymbol{\tau})}{\partial{x_3}}n_{3}\,,\label{eq:IGA-CNeumannK33}\\
&\tilde{\textbf{K}}_{23}(\boldsymbol{\tau})=\mathbb{C}_{23}\cfrac{\partial{R_{\textbf{i},\textbf{p}}}(\boldsymbol{\tau})}{\partial{x_3}}n_{2}+\mathbb{C}_{44}\cfrac{\partial{R_{\textbf{i},\textbf{p}}}(\boldsymbol{\tau})}{\partial{x_2}}n_{3}\,,\label{eq:IGA-CNeumannK23}\\
&\tilde{\textbf{K}}_{13}(\boldsymbol{\tau})=\mathbb{C}_{13}\cfrac{\partial{R_{\textbf{i},\textbf{p}}}(\boldsymbol{\tau})}{\partial{x_3}}n_{1}+\mathbb{C}_{55}\cfrac{\partial{R_{\textbf{i},\textbf{p}}}(\boldsymbol{\tau})}{\partial{x_1}}n_{3}\,,\label{eq:IGA-CNeumannK31}\\
&\tilde{\textbf{K}}_{12}(\boldsymbol{\tau})=\mathbb{C}_{12}\cfrac{\partial{R_{\textbf{i},\textbf{p}}}(\boldsymbol{\tau})}{\partial{x_2}}n_{1}+\mathbb{C}_{66}\cfrac{\partial{R_{\textbf{i},\textbf{p}}}(\boldsymbol{\tau})}{\partial{x_1}}n_{2}\,.\label{eq:IGA-CNeumannK12}
\end{align}
\end{subequations}
As we can see from equations~\eqref{eq:IGA-CNeumann}, Neumann boundary conditions are
directly imposed as strong equations at the collocation points belonging to the boundary surface (see,~\cite{Auricchio2012a,DeLorenzis2015}), with the usual physical meaning of prescribed boundary traction.

\section{An IGA collocation approach to model 3D composite plates}\label{sec:IGA-Cstrategies}
In this section we describe our IGA 3D collocation strategy to model composite plates.
The proposed method, known as single element approach, relies on a homogenization technique combined with a post-processing approach based on the imposition of equilibrium equations in strong form.	

\subsection{Single-element approach}\label{subsec:SEA}     
The single-element approach considers the plate discretized by a single element through the thickness, which strongly reduces the number of degrees of freedom with respect to layerwise methods. The material matrix is therefore homogenized to account for the presence of the layers as Figure~\ref{subfig-2:HSE} clearly describes. 
\begin{remark}
Considering a single-element homogenized approach is effective only for through-the-thickness symmetric layer distributions, 	as for non-symmetric ply stacking sequences the plate middle plane is not balanced. In the case of non-symmetric layer distributions this technique is still applicable when the stacking sequence can be split into two symmetric piles, using one element per homogenized stack with a $C^0$ interface.
\end{remark}
This method provides accurate results only in terms of displacements and in-plane stress components and, in order to recover a proper out-of-plane stress state, following~\cite{Dufour2018}, we propose to couple it with a post-processing technique. 
To characterize the variation of the material properties from layer to layer, we homogenize the constitutive behavior to create an equivalent single-layer laminate, referring to~\cite{Sun1988}, where explicit expressions for the effective elastic constants of the equivalent laminate are given as\newline
\begin{subequations}
\begin{align}
&\overline{\mathbb{C}}_{11}=\sum_{k=1}^{N}\overline{t}_{k}\mathbb{C}_{11}^{(k)}+\sum_{k=2}^{N}(\mathbb{C}_{13}^{(k)}-\overline{\mathbb{C}}_{13})\overline{t}_{k}\frac{(\mathbb{C}_{13}^{(1)}-\mathbb{C}_{13}^{(k)})}{\mathbb{C}_{33}^{(k)}}\label{eq:aveC1}\\
&\overline{\mathbb{C}}_{12}=\sum_{k=1}^{N}\overline{t}_{k}\mathbb{C}_{12}^{(k)}+\sum_{k=2}^{N}(\mathbb{C}_{13}^{(k)}-\overline{\mathbb{C}}_{13})\overline{t}_{k}\frac{(\mathbb{C}_{23}^{(1)}-\mathbb{C}_{23}^{(k)})}{\mathbb{C}_{33}^{(k)}}\label{eq:aveC2}\\
&\overline{\mathbb{C}}_{13}=\sum_{k=1}^{N}\overline{t}_{k}\mathbb{C}_{13}^{(k)}+\sum_{k=2}^{N}(\mathbb{C}_{33}^{(k)}-\overline{\mathbb{C}}_{33})\overline{t}_{k}\frac{(\mathbb{C}_{13}^{(1)}-\mathbb{C}_{13}^{(k)})}{\mathbb{C}_{33}^{(k)}}\label{eq:aveC3}\\
&\overline{\mathbb{C}}_{22}=\sum_{k=1}^{N}\overline{t}_{k}\mathbb{C}_{22}^{(k)}+\sum_{k=2}^{N}(\mathbb{C}_{23}^{(k)}-\overline{\mathbb{C}}_{23})\overline{t}_{k}\frac{(\mathbb{C}_{23}^{(1)}-\mathbb{C}_{23}^{(k)})}{\mathbb{C}_{33}^{(k)}}\label{eq:aveC4}\\
&\overline{\mathbb{C}}_{23}=\sum_{k=1}^{N}\overline{t}_{k}\mathbb{C}_{23}^{(k)}+\sum_{k=2}^{N}(\mathbb{C}_{33}^{(k)}-\overline{\mathbb{C}}_{33})\overline{t}_{k}\frac{(\mathbb{C}_{23}^{(1)}-\mathbb{C}_{23}^{(k)})}{\mathbb{C}_{33}^{(k)}}\label{eq:aveC5}\\
&\overline{\mathbb{C}}_{33}=\frac{1}{\bigg(\sum_{k=1}^{N}\cfrac{\overline{t}_{k}}{\mathbb{C}_{33}^{(k)}}\bigg)}\label{eq:aveC6}\\
&\overline{\mathbb{C}}_{44}=\frac{\bigg(\sum_{k=1}^{N}\cfrac{\overline{t}_{k}\mathbb{C}_{44}^{(k)}}{\Delta_{k}}\bigg)}{\Delta},\hspace{0.5cm}\Delta=\bigg(\sum_{k=1}^{N}\frac{\overline{t}_{k}\mathbb{C}_{44}^{(k)}}{\Delta_{k}}\bigg)\bigg(\sum_{k=1}^{N}\frac{\overline{t}_{k}\mathbb{C}_{55}^{(k)}}{\Delta_{k}}\bigg)\label{eq:aveC7}\\
&\overline{\mathbb{C}}_{55}=\frac{\bigg(\sum_{k=1}^{N}\cfrac{\overline{t}_{k}\mathbb{C}_{55}^{(k)}}{\Delta_{k}}\bigg)}{\Delta},\hspace{0.5cm}\Delta_{k}=\mathbb{C}_{44}^{k}\mathbb{C}_{55}^{k}\label{eq:aveC8}\\
&\overline{\mathbb{C}}_{66}=\sum_{k=1}^{N}\overline{t}_{k}\mathbb{C}_{66}^{(k)}\label{eq:aveC9}
\end{align}
\end{subequations}
where $\mathbb{C}_{ij}^{(k)}$ represents the $ij$-th component of the fourth order elasticity tensor in Voigt notation \eqref{eq:CIVeltensor} for the $k$-th layer and $\overline{t}_{k}=\cfrac{t_{k}}{h}$ stands for the volume fraction of the $k$-th lamina, $h$ being the total thickness and $t_{k}$ the $k$-th thickness.

\subsubsection{Post-processing step: Reconstruction from Equilibrium}\label{subsec:post-processing}
As interlaminar delamination and other fracture processes rely mostly on out-of-plane components, a proper through-the-thickness stress description is required. In order to recover a more accurate stress state, we perform a post-processing step based on the equilibrium equations, following~\cite{Dufour2018}, relying on the higher regularity granted by IGA shape functions. This procedure, which takes its roots in~\cite{deMiranda2002,Engblom1985,Pryor1971,Ubertini2004}, has already been proven to be successful for IGA-Galerkin. Inside the plate the stresses should satisfy the equilibrium equation~\eqref{eq:equilibriumint} that can be expanded as
\begin{subequations}
\begin{align}
&\sigma_{11,1}+\sigma_{12,2}+\sigma_{13,3}=-b_{1}\,,\label{eq:equilibriumeng1}\\
&\sigma_{12,1}+\sigma_{22,2}+\sigma_{23,3}=-b_{2}\,,\label{eq:equilibriumeng2}\\
&\sigma_{13,1}+\sigma_{23,2}+\sigma_{33,3}=-b_{3}\,.\label{eq:equilibriumeng3}
\end{align}\label{eq:equilibriumeng}%
\end{subequations}
Assuming the in-plane stress components to well approximate the laminate behaviour, as it will be shown in Section~\ref{sec:results}, we can integrate equation~\ref{eq:equilibriumeng1} and~\ref{eq:equilibriumeng2} along the thickness, recovering the out-of-plane shear stresses as
  \begin{subequations}
	\begin{align}
	\label{eq:ppsigma13}
	\sigma_{13}(X_3) &= -\int^{X_3}_{\bar{X_3}}(\sigma_{11,1}(\zeta) + \sigma_{12,2}(\zeta)+b_1(\zeta))\rmd \zeta + \sigma_{13}(\bar{X_3})\,,\\
	\label{eq:ppsigma23}
	\sigma_{23}(X_3) &= -\int^{X_3}_{\bar{X_3}}(\sigma_{12,1}(\zeta) + \sigma_{22,2}(\zeta)+b_2(\zeta))\rmd \zeta + \sigma_{23}(\bar{X_3})\,,
	\end{align}
\end{subequations}
where $\zeta$ represents the coordinate along the thickness direction.\newline
Finally we can insert equations~\eqref{eq:ppsigma13} and~\eqref{eq:ppsigma23} into~\eqref{eq:equilibriumeng3}, recovering the $\sigma_{33}$ component as
\begin{equation}
\begin{aligned}
\label{eq:ppsigma33}
\sigma_{33}(X_3) &= -\int^{X_3}_{\bar{X_3}}(\sigma_{13,1}(\zeta) + \sigma_{23,2}(\zeta)+b_3(\zeta))\rmd \zeta + \sigma_{33}(\bar{X_3})\,.
\end{aligned}
\end{equation}
Following~\cite{Dufour2018}, the integral constants are chosen to fulfil the boundary conditions at the top or bottom surfaces.\newline
Recalling that
\begin{equation}
\begin{aligned}
&\sigma_{ij,k}=\overline{\mathbb{C}}_{ijmn}\frac{u_{m,nk}+u_{n,mk}}{2}\,,\label{eq:equilibriumsigmau}
\end{aligned}
\end{equation}
where the homogenized elasticity tensor $\overline{\mathbb{C}}$ is constant, it is clear the necessity of a highly regular displacement solution in order to recover a proper stress state. Such a condition can be easily achieved using isogeometric collocation, due to  the possibility to benefit from the high regularity of B-Splines or NURBS. We also remark that the proposed method strongly relies on the possibility to obtain an accurate description (with a relatively coarse mesh) of the in-plane stress state.

\section{Numerical tests}\label{sec:results}
In this section, to assess whether the proposed method can effectively reproduce composite plates behaviour, we consider a classical benchmark problem~\cite{Pagano1970} and we address different aspects such as the effectiveness of the proposed post-processing step, the method sensitivity to parameters of interest (i.e., number of layers and length-to-thickness ratio), and its convergence. 

\subsection{Reference solution: the Pagano layered plate}\label{subsec:Pagano}
\begin{figure}[!htbp]
	\centering
	\includegraphics[width=.6\textwidth]{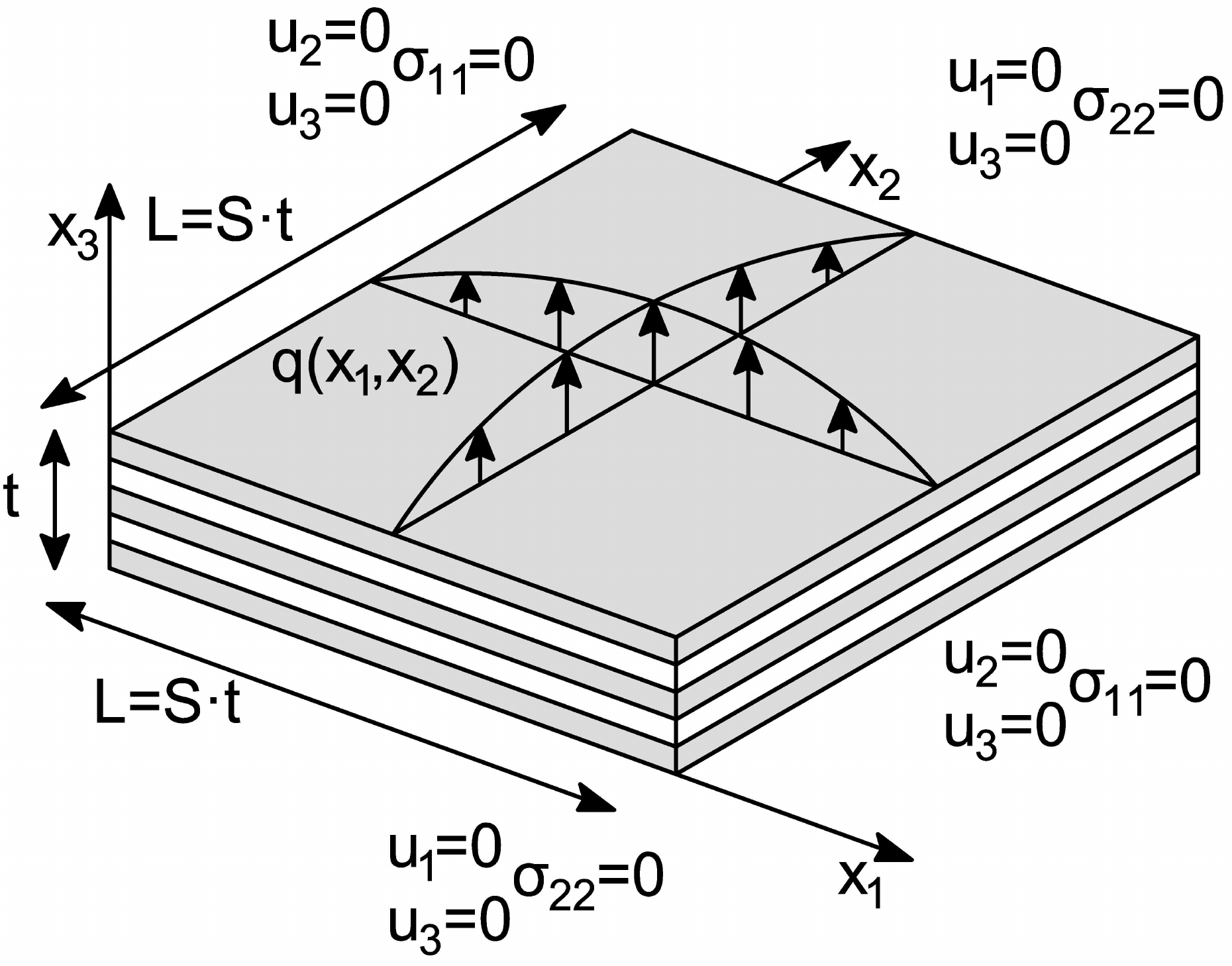}
	\caption{Pagano's test case~\cite{Pagano1970}. Problem geometry and boundary conditions.}
	\label{fig:testproblem}
\end{figure}

A square laminated composite plate of total thickness $t$ made of $N$ orthotropic layers is considered. This structure is simply supported and a
normal sinusoidal traction is applied on the upper surface, while the lower surface is traction-free, as shown in Figure~\ref{fig:testproblem}. 
In the proposed numerical tests we consider different numbers of layers, namely 3, 11 and 33.
The thickness of every single layer is set to 1 mm, and the edge length, $L$, is chosen to be $S$ times larger than the total thickness $t$ of the laminate.
Different choices of length-to-thickness ratio are considered (i.e., 20, 30, 40, and 50) which allow to draw interesting considerations about the laminate behaviour in the proposed convergence tests.
For all examples we consider the same loading conditions proposed by Pagano, i.e., a double sinus with periodicity equal to twice the length of the plate. 
As depicted in Figure~\ref{fig:testproblem} the laminated plate is composed of layers organized in an alternated distribution of orthotropic plies (i.e., a 0\textdegree/90\textdegree~stacking sequence in our case).
Layer material parameters considered in the numerical tests are summarized in Table~\ref{tab:matproperties} for 0\textdegree-oriented plies.
\begin{table}[!htbp]
	\caption{Material properties for 0\textdegree-oriented layers employed in the numerical tests.} 
	\begin{center}
		\small{
		\begin{tabular}{?c c c c c c c c c?}
			\thickhline 
			$E_{1}$ & $E_{2}$& $E_{3}$ & $G_{23}$ & $G_{13}$ & $G_{12}$ & $\nu_{23}$ & $\nu_{13}$ & $\nu_{12}$\Tstrut\Bstrut\\\hline 
			[GPa] & [GPa] & [GPa] & [GPa] & [GPa] & [GPa]& [-] & [-] & [-]\Tstrut\Bstrut\\\thickhline 
			25000 & 1000 & 1000 & 200 & 500 & 500 & 0.25 & 0.25& 0.25\Tstrut\Bstrut\\\thickhline  
		\end{tabular}}
	\end{center}\label{tab:matproperties}
\end{table}
\newpage
The Neumann boundary conditions on the plate surfaces $x_3=\pm\dfrac{t}{2}$ are
\begin{equation}
\begin{aligned}
&\sigma_{33}(x_1,x_2,-\dfrac{t}{2})=\sigma_{13}(x_1,x_2,\pm\dfrac{t}{2})=\sigma_{23}(x_1,x_2,\pm\dfrac{t}{2})=0\,,\\
&\sigma_{33}(x_1,x_2,\dfrac{t}{2})=\sigma_0\sin(\dfrac{\pi x_1}{St})\sin(\dfrac{\pi x_2}{St})\,,\label{eq:NBCs}
\end{aligned}
\end{equation}
where $\sigma_0 =$ 1 MPa.\newline
The simple support edge conditions are taken as
\begin{equation}
\begin{aligned}
\bullet\;\sigma_{11}=0\;\text{and}\;u_2=u_3=0\;\text{at}\;x_1=0\;\text{and}\;x_1=L\,,\\
\bullet\;\sigma_{22}=0\;\text{and}\;u_1=u_3=0\;\text{at}\;x_2=0\;\text{and}\;x_2=L\,.\label{eq:DBCs}
\end{aligned}
\end{equation}
All results are then expressed in terms of the following normalized stress components
\begin{equation}
\begin{aligned}
&\overline{\sigma}_{ij}=\dfrac{\sigma_{ij}}{\sigma_0S^2},\hspace{1cm}i,j=1,2\,,\\
&\overline{\sigma}_{i3}=\dfrac{\sigma_{i3}}{\sigma_0S},\hspace{1.15cm}i=1,2\,,\\
&\overline{\sigma}_{33}=\dfrac{\sigma_{33}}{\sigma_0}\,.\\\label{eq:normalizedresults}
\end{aligned}
\end{equation}

\subsection{Post-processed out-of-plane stresses}\label{subsec:ppeffects}
\begin{figure}[!htbp]
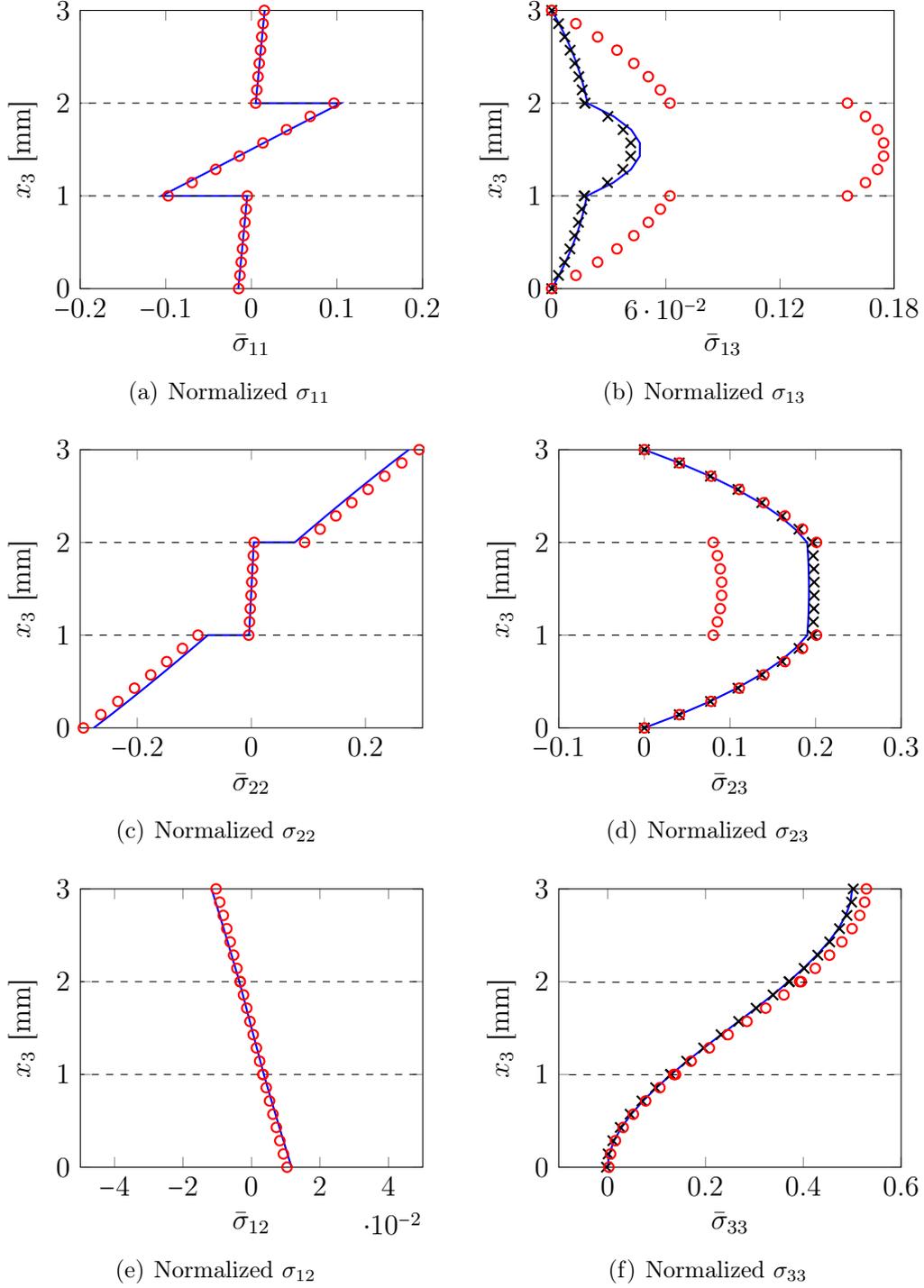

	\centering
	\subfigure[Normalized $\sigma_{11}$\label{subfig-1:test3}]{\ifrecompiletikz\tikzsetnextfilename{fig_02_a}\tikzexternalenable\begin{tikzpicture}
\begin{axis}[stress plot style,
  xmin=-0.2, xmax=0.2,
  ymin=0, ymax=3,
  xlabel={$\bar{\sigma}_{11}$}]
  
\addplot[mark=none, blue,thick] table[x index=2,y index=0, col sep=comma] {srcTikz/S11_lay3/Sig11_1_1.dat};
\addplot[red,mark=o, mark size = 2, only marks,thick] table[x index=1,y index=0, col sep=comma] {srcTikz/S11_lay3/Sig11_1_1.dat};

\foreach \lin in {1,...,2}{
	\addplot[dashed,mark=none, black, samples=2] {\lin};
}
\end{axis}
\end{tikzpicture}\tikzexternaldisable\else\includegraphics{fig_02_a.pdf}\fi}\hfill
	\subfigure[Normalized $\sigma_{13}$\label{subfig-2:test3}]{\ifrecompiletikz\tikzsetnextfilename{fig_02_b}\tikzexternalenable\begin{tikzpicture}
\begin{axis}[stress plot style,
  xmin=0, xmax=0.18,
  ymin=0, ymax=3, xtick={0,0.06,0.12,0.18},
  xlabel={$\bar{\sigma}_{13}$}]
  
\addplot[mark=none, blue,thick] table[x index=2,y index=0, col sep=comma] {srcTikz/S13_lay3/Sig13_1_1.dat};
\addplot[black,mark=x, mark size = 3, only marks,thick] table[x index=1,y index=0, col sep=comma] {srcTikz/S13_lay3/Sig13_1_1.dat};
\addplot[red,mark=o, mark size = 2, only marks,thick] table[x index=3,y index=0, col sep=comma] {srcTikz/S13_lay3/Sig13_1_1.dat};

\foreach \lin in {1,...,2}{
	\addplot[dashed,mark=none, black, samples=2] {\lin};
}
\end{axis}
\end{tikzpicture}\tikzexternaldisable\else\includegraphics{fig_02_b.pdf}\fi}\\
	\subfigure[Normalized $\sigma_{22}$\label{subfig-3:test3}]{\ifrecompiletikz\tikzsetnextfilename{fig_02_c}\tikzexternalenable\begin{tikzpicture}
\begin{axis}[stress plot style,
  xmin=-0.3, xmax=0.3,
  ymin=0, ymax=3,
  xlabel={$\bar{\sigma}_{22}$}]
  
\addplot[mark=none, blue,thick] table[x index=2,y index=0, col sep=comma] {srcTikz/S22_lay3/Sig22_1_1.dat};
\addplot[red,mark=o, mark size = 2, only marks,thick] table[x index=1,y index=0, col sep=comma] {srcTikz/S22_lay3/Sig22_1_1.dat};

\foreach \lin in {1,...,2}{
	\addplot[dashed,mark=none, black, samples=2] {\lin};
}
\end{axis}
\end{tikzpicture}\tikzexternaldisable\else\includegraphics{fig_02_c.pdf}\fi}\hfill
	\subfigure[Normalized $\sigma_{23}$\label{subfig-4:test3}]{\ifrecompiletikz\tikzsetnextfilename{fig_02_d}\tikzexternalenable\begin{tikzpicture}
\begin{axis}[stress plot style,
  xmin=-0.1, xmax=0.3,
  ymin=0, ymax=3,
  xlabel={$\bar{\sigma}_{23}$}]
  
\addplot[mark=none, blue,thick] table[x index=2,y index=0, col sep=comma] {srcTikz/S23_lay3/Sig23_1_1.dat};
\addplot[black,mark=x, mark size = 3, only marks,thick] table[x index=1,y index=0, col sep=comma] {srcTikz/S23_lay3/Sig23_1_1.dat};
\addplot[red,mark=o, mark size = 2, only marks,thick] table[x index=3,y index=0, col sep=comma] {srcTikz/S23_lay3/Sig23_1_1.dat};

\foreach \lin in {1,...,2}{
	\addplot[dashed,mark=none, black, samples=2] {\lin};
}
\end{axis}
\end{tikzpicture}\tikzexternaldisable\else\includegraphics{fig_02_d.pdf}\fi}\\
	\subfigure[Normalized $\sigma_{12}$\label{subfig-5:test3}]{\ifrecompiletikz\tikzsetnextfilename{fig_02_e}\tikzexternalenable\begin{tikzpicture}
\begin{axis}[stress plot style,
  xmin=-0.05, xmax=0.05,
  ymin=0, ymax=3,
  xlabel={$\bar{\sigma}_{12}$}]
  
\addplot[mark=none, blue,thick] table[x index=2,y index=0, col sep=comma] {srcTikz/S12_lay3/Sig12_1_1.dat};
\addplot[red,mark=o, mark size = 2, only marks,thick] table[x index=1,y index=0, col sep=comma] {srcTikz/S12_lay3/Sig12_1_1.dat};

\foreach \lin in {1,...,2}{
	\addplot[dashed,mark=none, black, samples=2] {\lin};
}
\end{axis}
\end{tikzpicture}\tikzexternaldisable\else\includegraphics{fig_02_e.pdf}\fi}\hfill
	\subfigure[Normalized $\sigma_{33}$\label{subfig-6:test3}]{\ifrecompiletikz\tikzsetnextfilename{fig_02_f}\tikzexternalenable\begin{tikzpicture}
\begin{axis}[stress plot style,
  xmin=-0.1, xmax=0.6,
  ymin=0, ymax=3,
  xlabel={$\bar{\sigma}_{33}$}]
 
\addplot[mark=none, blue,thick] table[x index=2,y index=0, col sep=comma] {srcTikz/S33_lay3/Sig33_1_1.dat};
\addplot[black,mark=x, mark size = 3, only marks,thick] table[x index=1,y index=0, col sep=comma] {srcTikz/S33_lay3/Sig33_1_1.dat};
\addplot[red,mark=o, mark size = 2, only marks,thick] table[x index=3,y index=0, col sep=comma] {srcTikz/S33_lay3/Sig33_1_1.dat};
\end{axis}

\draw [dashed] (0.15,2.72) -- (5,2.72);
\draw [dashed] (0.15,1.36) -- (5,1.36);
\end{tikzpicture}\tikzexternaldisable\else\includegraphics{fig_02_f.pdf}\fi}
	\caption{Through-the-thickness stress solutions for the 3D Pagano problem~\cite{Pagano1970} evaluated at $x_1=x_2=0.25L$. Case: plate with 3 layers and length-to-thickness ratio $S=20$ (\textcolor{blue}{$\boldsymbol{\leftrightline}$} Pagano's solution, \protect\markerfive\hspace{0.2cm}homogenized single-element approach solution (without post-processing), $\boldsymbol{\times}$ post-processed solution).}
	\label{fig:example_sig3}
\end{figure}

\begin{figure}
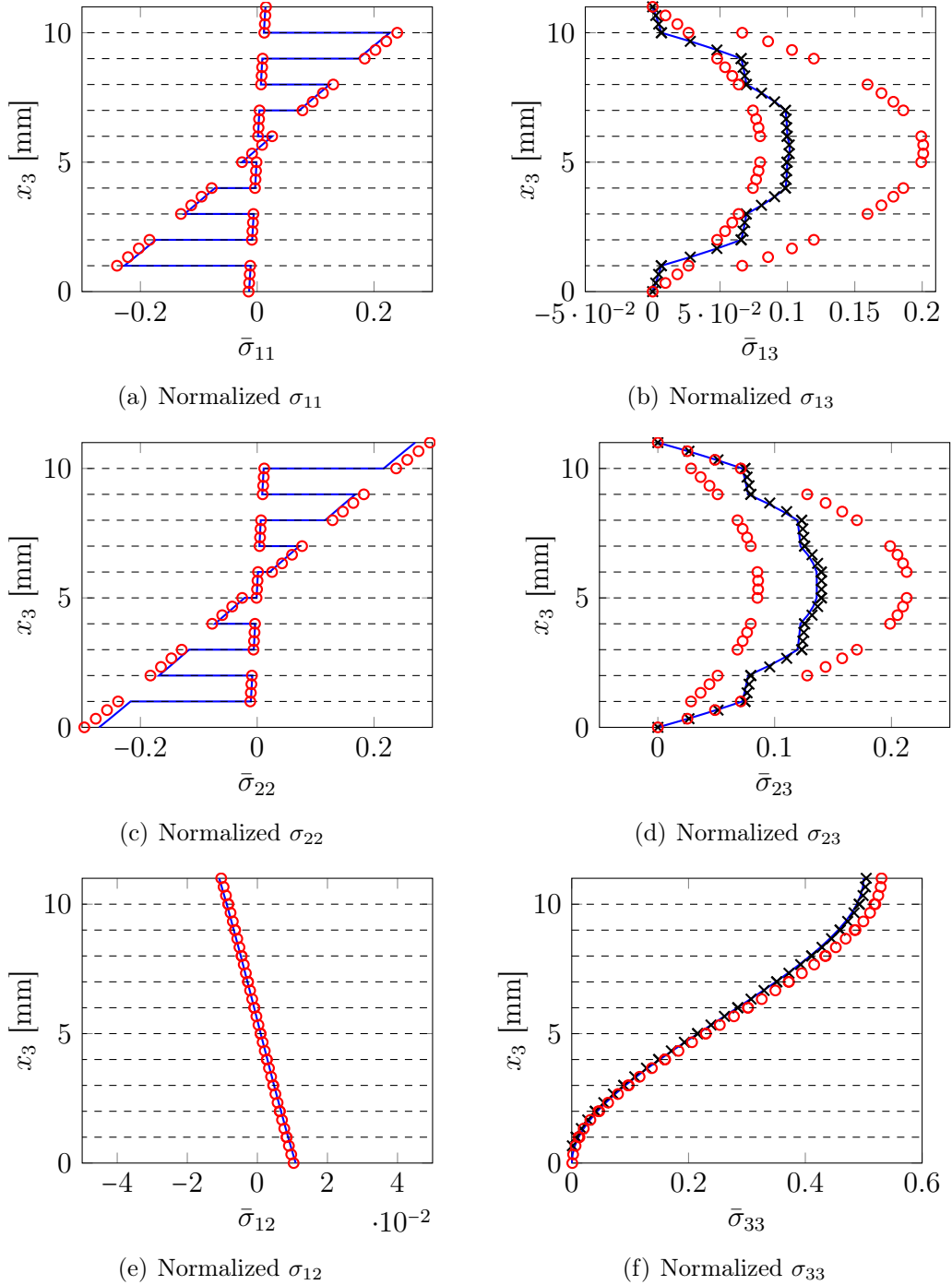

	\centering
	\subfigure[Normalized $\sigma_{11}$\label{subfig-1:test11}]{\ifrecompiletikz\tikzsetnextfilename{fig_03_a}\tikzexternalenable\begin{tikzpicture}
\begin{axis}[stress plot style,
  xmin=-0.3, xmax=0.3,
  ymin=0, ymax=11,
  xlabel={$\bar{\sigma}_{11}$}]
  
\addplot[mark=none, blue,thick] table[x index=2,y index=0, col sep=comma] {srcTikz/S11_lay11/Sig11_1_1.dat};
\addplot[red,mark=o, mark size = 2, only marks,thick] table[x index=1,y index=0, col sep=comma,each nth point={1}] {srcTikz/S11_lay11/Sig11_1_1.dat};

\foreach \lin in {1,...,10}{
	\addplot[dashed,mark=none, black, samples=2] {\lin};
}
\end{axis}
\end{tikzpicture}\tikzexternaldisable\else\includegraphics{fig_03_a.pdf}\fi}\hfill
	\subfigure[Normalized $\sigma_{13}$\label{subfig-2:test11}]{\ifrecompiletikz\tikzsetnextfilename{fig_03_b}\tikzexternalenable\begin{tikzpicture}
\begin{axis}[stress plot style,
  xmin=-0.05, xmax=0.21,
  ymin=0, ymax=11,
  xlabel={$\bar{\sigma}_{13}$}]
  
\addplot[mark=none, blue,thick] table[x index=2,y index=0, col sep=comma] {srcTikz/S13_lay11/Sig13_1_1.dat};
\addplot[black,mark=x, mark size = 3, only marks,thick] table[x index=1,y index=0, col sep=comma,each nth point={1}] {srcTikz/S13_lay11/Sig13_1_1.dat};
\addplot[red,mark=o, mark size = 2, only marks,thick] table[x index=3,y index=0, col sep=comma,each nth point={1}] {srcTikz/S13_lay11/Sig13_1_1.dat};

\foreach \lin in {1,...,10}{
	\addplot[dashed,mark=none, black, samples=2] {\lin};
}
\end{axis}
\end{tikzpicture}\tikzexternaldisable\else\includegraphics{fig_03_b.pdf}\fi}\\
	\subfigure[Normalized $\sigma_{22}$\label{subfig-3:test11}]{\ifrecompiletikz\tikzsetnextfilename{fig_03_c}\tikzexternalenable\begin{tikzpicture}
\begin{axis}[stress plot style,
  xmin=-0.3, xmax=0.3,
  ymin=0, ymax=11,
  xlabel={$\bar{\sigma}_{22}$}]
  
\addplot[mark=none, blue,thick] table[x index=2,y index=0, col sep=comma] {srcTikz/S22_lay11/Sig22_1_1.dat};
\addplot[red,mark=o, mark size = 2, only marks,thick] table[x index=1,y index=0, col sep=comma,each nth point={1}] {srcTikz/S22_lay11/Sig22_1_1.dat};

\foreach \lin in {1,...,10}{
	\addplot[dashed,mark=none, black, samples=2] {\lin};
}
\end{axis}
\end{tikzpicture}\tikzexternaldisable\else\includegraphics{fig_03_c.pdf}\fi}\hfill
	\subfigure[Normalized $\sigma_{23}$\label{subfig-4:test11}]{\ifrecompiletikz\tikzsetnextfilename{fig_03_d}\tikzexternalenable\begin{tikzpicture}
\begin{axis}[stress plot style,
  xmin=-0.05, xmax=0.25,
  ymin=0, ymax=11,
  xlabel={$\bar{\sigma}_{23}$}]
  
\addplot[mark=none, blue,thick] table[x index=2,y index=0, col sep=comma] {srcTikz/S23_lay11/Sig23_1_1.dat};
\addplot[black,mark=x, mark size = 3, only marks,thick] table[x index=1,y index=0, col sep=comma,each nth point={1}] {srcTikz/S23_lay11/Sig23_1_1.dat};
\addplot[red,mark=o, mark size = 2, only marks,thick] table[x index=3,y index=0, col sep=comma,each nth point={1}] {srcTikz/S23_lay11/Sig23_1_1.dat};

\foreach \lin in {1,...,10}{
	\addplot[dashed,mark=none, black, samples=2] {\lin};
}
\end{axis}
\end{tikzpicture}\tikzexternaldisable\else\includegraphics{fig_03_d.pdf}\fi}\\
	\subfigure[Normalized $\sigma_{12}$\label{subfig-5:test11}]{\ifrecompiletikz\tikzsetnextfilename{fig_03_e}\tikzexternalenable\begin{tikzpicture}
\begin{axis}[stress plot style,
  xmin=-0.05, xmax=0.05,
  ymin=0, ymax=11,
  xlabel={$\bar{\sigma}_{12}$}]
  
\addplot[mark=none, blue,thick] table[x index=2,y index=0, col sep=comma] {srcTikz/S12_lay11/Sig12_1_1.dat};
\addplot[red,mark=o, mark size = 2, only marks,thick] table[x index=1,y index=0, col sep=comma,each nth point={1}] {srcTikz/S12_lay11/Sig12_1_1.dat};

\foreach \lin in {1,...,10}{
	\addplot[dashed,mark=none, black, samples=2] {\lin};
}
\end{axis}
\end{tikzpicture}\tikzexternaldisable\else\includegraphics{fig_03_e.pdf}\fi}\hfill
	\subfigure[Normalized $\sigma_{33}$\label{subfig-6:test11}]{\ifrecompiletikz\tikzsetnextfilename{fig_03_f}\tikzexternalenable\begin{tikzpicture}
\begin{axis}[stress plot style,
  xmin=0, xmax=0.6,
  ymin=0, ymax=11,
  xlabel={$\bar{\sigma}_{33}$}]
\addplot[mark=none, blue,thick] table[x index=2,y index=0, col sep=comma] {srcTikz/S33_lay11/Sig33_1_1.dat};
\addplot[black,mark=x, mark size = 3, only marks,thick] table[x index=1,y index=0, col sep=comma,each nth point={1}] {srcTikz/S33_lay11/Sig33_1_1.dat};
\addplot[red,mark=o, mark size = 2, only marks,thick] table[x index=3,y index=0, col sep=comma,each nth point={1}] {srcTikz/S33_lay11/Sig33_1_1.dat};
\end{axis}

\draw [dashed] (0,0.372) -- (5,0.372);
\draw [dashed] (0,0.744) -- (5,0.744);
\draw [dashed] (0,1.116) -- (5,1.116);
\draw [dashed] (0,1.488) -- (5,1.488);
\draw [dashed] (0,1.86) -- (5,1.86);
\draw [dashed] (0,2.232) -- (5,2.232);
\draw [dashed] (0,2.604) -- (5,2.604);
\draw [dashed] (0,2.976) -- (5,2.976);
\draw [dashed] (0,3.348) -- (5,3.348);
\draw [dashed] (0,3.72) -- (5,3.72);
\end{tikzpicture}\tikzexternaldisable\else\includegraphics{fig_03_f.pdf}\fi}
	\caption{Through-the-thickness stress solutions for the 3D Pagano problem~\cite{Pagano1970} evaluated at $x_1=x_2=0.25L$. Case: plate with 11 layers and length-to-thickness ratio $S=20$ (\textcolor{blue}{$\boldsymbol{\leftrightline}$} Pagano's solution, \protect\markerfive\hspace{0.2cm}homogenized single-element approach solution (without post-processing), $\boldsymbol{\times}$ post-processed solution).}
	\label{fig:example_sig11}
\end{figure}

In this section, we comment the results obtained using the proposed IGA-collocation approach, as compared with Pagano's analytical solution~\cite{Pagano1970}. To give an idea of the improvement granted by the post-processing of out-of-plane stress components, in Figures~\ref{fig:example_sig3} and~\ref{fig:example_sig11} we compare the reference solution with non-treated and post-processed results
for the cases with 3 and 11 layers, considering a length to thickness ratio $S=20$.
All numerical simulations are carried out using an in-plane degree of approximation $p=q=6$ and 10 collocation points for each in-plane parametric direction, while we use an approximation degree $r=4$ and one element through the thickness (i.e., $r+1$ collocation points).
The sampling point where we show results is located at $x_1=x_2=0.25L$.
For both considered cases the in-plane stresses show a good behaviour, as expected, while the out-of-plane stress components, without a post-processing treatment, are erroneusly discontinuous. 
The proposed results clearly show the improvement granted by the post-process of out-of-plane components.

\begin{figure}[!htbp]
		\centering
    \ifrecompiletikz\tikzsetnextfilename{fig_04}\tikzexternalenable\begin{tikzpicture}
    \node at (1.1, -11.5) {\small $x_1=0$};
    \node at (3.7, -11.5) {\small $x_1=L/4$};
    \node at (6.3, -11.5) {\small $x_1=L/2$};
    \node at (8.9, -11.5) {\small $x_1=3\,L/4$};
    \node at (11.5, -11.5) {\small $x_1=L$};

    \node [rotate=90] at (-1, 0.8) {\small $x_2=L$};
    \node [rotate=90] at (-1, -1.8) {\small $x_2=3\,L/4$};
    \node [rotate=90] at (-1, -4.4) {\small $x_2=L/2$};
    \node [rotate=90] at (-1, -7) {\small $x_2=L/4$};
    \node [rotate=90] at (-1, -9.6) {\small $x_2=0$};

    \begin{groupplot}[group style={group size=5 by 5,horizontal sep=.5cm},ytick={0,5,11},width=.27\textwidth,xmin=-0.25, xmax=0.25,enlarge y limits = false, axis x line=bottom,axis y line=left]
		\def\myPlots{}
			\pgfplotsforeachungrouped \y in {4,3,2,1,0}
				{
				\pgfplotsforeachungrouped \x in {0,1,2,3,4}
					{
						\eappto\myPlots{%
						\noexpand\nextgroupplot
						\noexpand\addplot+[red,thick, mark=none]
						table[x index=1, y index=0, col sep=comma] {srcTikz/S13_sampling/Sig13_\x_\y.dat};
						\noexpand\addplot+[only marks,blue, mark=x,each nth point={9}]
						table[x index=2, y index=0, col sep=comma] {srcTikz/S13_sampling/Sig13_\x_\y.dat};
					}
				}
		}
		\myPlots
	\end{groupplot}
\end{tikzpicture}\tikzexternaldisable\else\includegraphics{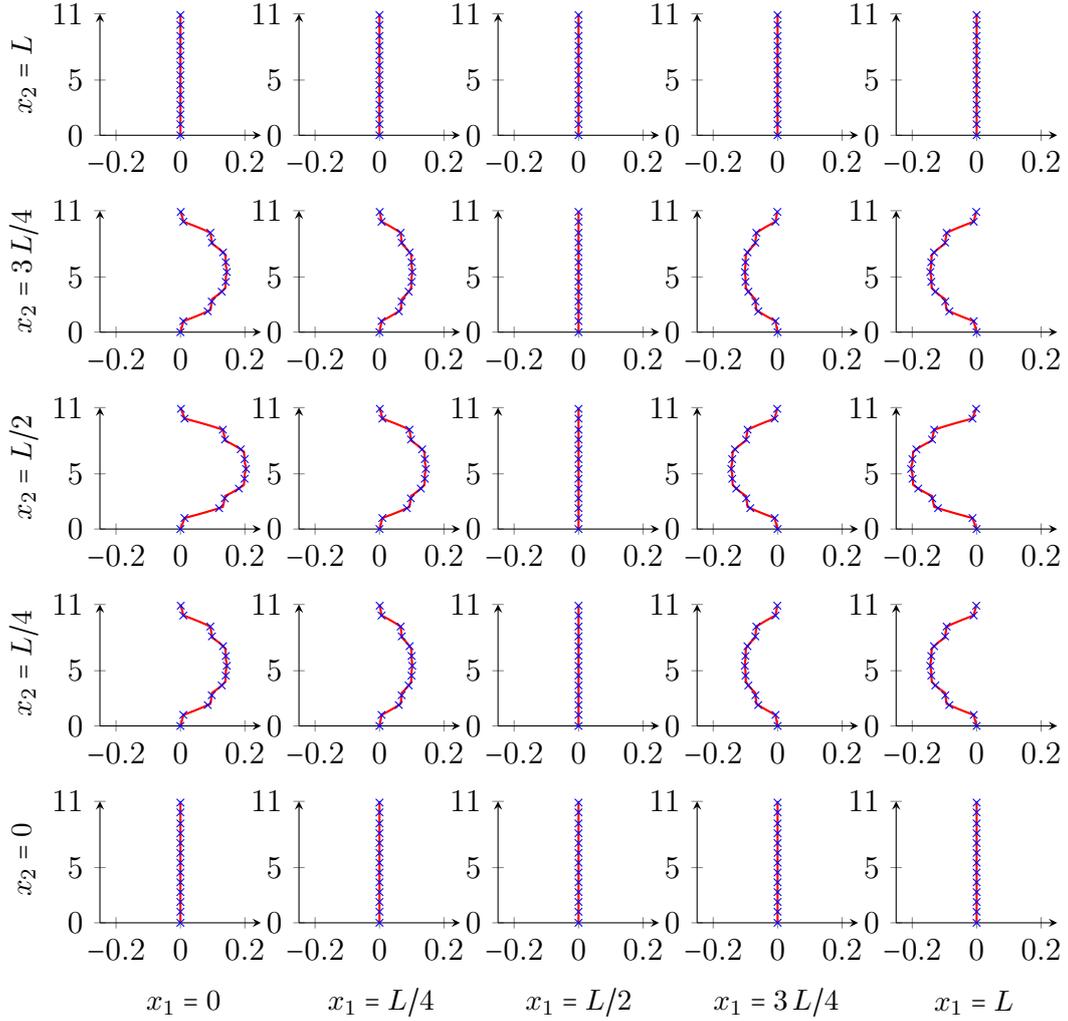}\fi
		\caption{Through-the-thickness $\bar{\sigma}_{13}$ profiles for several in plane sampling points.
		$L$ represents the total length of the plate, that for this case is  $L=220\,\text{mm}$ (being $L=S\,t$ with $t=11\,\text{mm}$ and $S=20$), while the number of layers is 11 (\textcolor{red}{$\boldsymbol{\leftrightline}$}post-processed solution, \textcolor{blue}{$\boldsymbol{\times}$} analytical solution~\cite{Pagano1970}).}
		\label{fig:samplingS13}
	\end{figure}

\begin{figure}[!htbp]
		\centering
    \ifrecompiletikz\tikzsetnextfilename{fig_05}\tikzexternalenable\begin{tikzpicture}
    \node at (1.1, -11.5) {\small $x_1=0$};
    \node at (3.7, -11.5) {\small $x_1=L/4$};
    \node at (6.3, -11.5) {\small $x_1=L/2$};
    \node at (8.9, -11.5) {\small $x_1=3\,L/4$};
    \node at (11.5, -11.5) {\small $x_1=L$};

    \node [rotate=90] at (-1, 0.8) {\small $x_2=L$};
    \node [rotate=90] at (-1, -1.8) {\small $x_2=3\,L/4$};
    \node [rotate=90] at (-1, -4.4) {\small $x_2=L/2$};
    \node [rotate=90] at (-1, -7) {\small $x_2=L/4$};
    \node [rotate=90] at (-1, -9.6) {\small $x_2=0$};

      \begin{groupplot}[group style={group size=5 by 5,horizontal sep=.5cm},ytick={0,5,11},width=.27\textwidth,xmin=-0.35, xmax=0.35, enlarge y limits = false, axis x line=bottom,axis y line=left]
		\def\myPlots{}
		\pgfplotsforeachungrouped \y in {4,3,2,1,0}
		{
			\pgfplotsforeachungrouped \x in {0,1,2,3,4}
			{
				\eappto\myPlots{%
					\noexpand\nextgroupplot
					\noexpand\addplot+[red,thick, mark=none]
					table[x index=1, y index=0, col sep=comma] {srcTikz/S23_sampling/Sig23_\x_\y.dat};
					\noexpand\addplot+[only marks,blue, mark=x,each nth point={9}]
					table[x index=2, y index=0, col sep=comma] {srcTikz/S23_sampling/Sig23_\x_\y.dat};
				}
			}
		}
		\myPlots
		\end{groupplot}
		\end{tikzpicture}\tikzexternaldisable\else\includegraphics{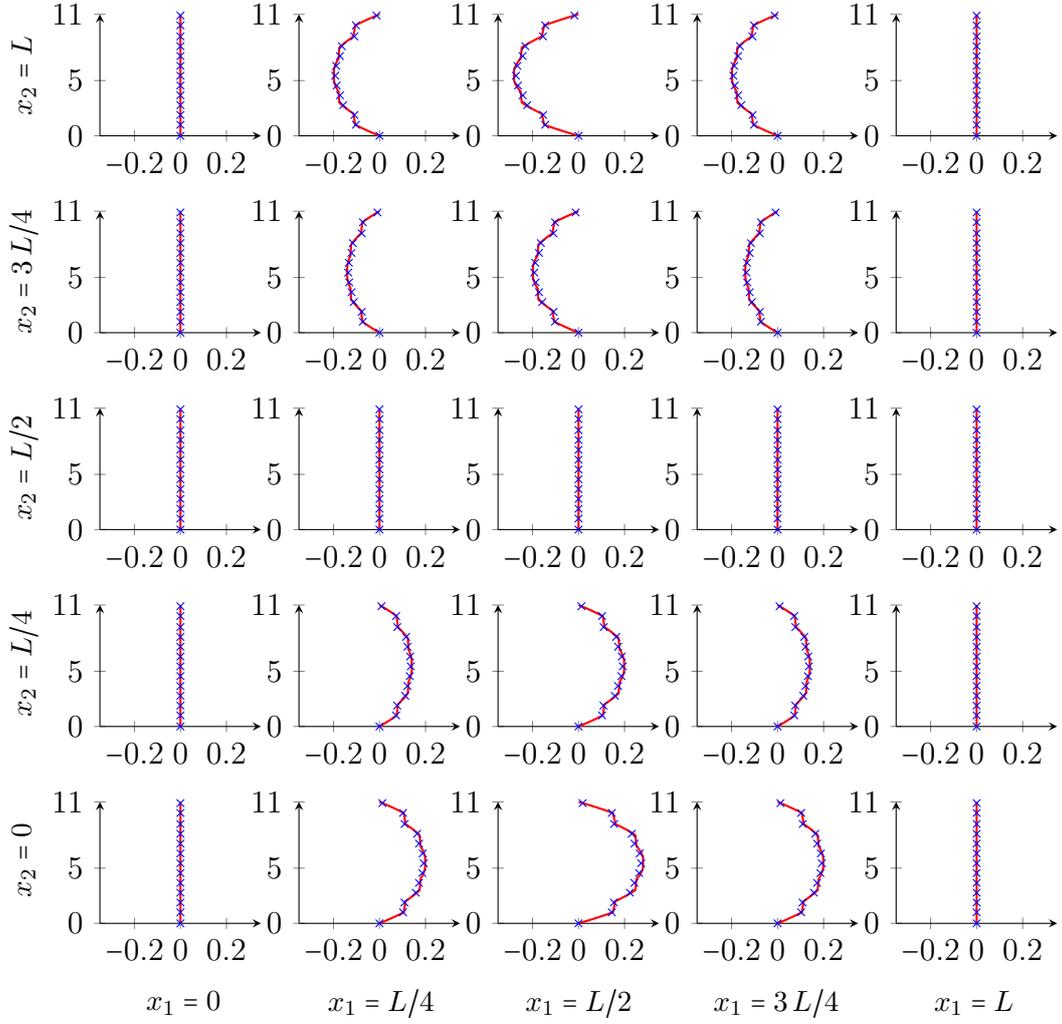}\fi
		\caption{Through-the-thickness $\bar{\sigma}_{23}$ profiles for several in plane sampling points.
		$L$ represents the total length of the plate, that for this case is  $L=220\,\text{mm}$ (being $L=S\,t$ with $t=11\,\text{mm}$ and $S=20$), while the number of layers is 11 (\textcolor{red}{$\boldsymbol{\leftrightline}$}post-processed solution, \textcolor{blue}{$\boldsymbol{\times}$} analytical solution~\cite{Pagano1970}).}
		\label{fig:samplingS23}
	\end{figure}

\begin{figure}[!htbp]
	\centering
  \ifrecompiletikz\tikzsetnextfilename{fig_06}\tikzexternalenable\begin{tikzpicture}
    \node at (1.1, -11.5) {\small $x_1=0$};
    \node at (3.7, -11.5) {\small $x_1=L/4$};
    \node at (6.3, -11.5) {\small $x_1=L/2$};
    \node at (8.9, -11.5) {\small $x_1=3\,L/4$};
    \node at (11.5, -11.5) {\small $x_1=L$};

    \node [rotate=90] at (-1, 0.8) {\small $x_2=L$};
    \node [rotate=90] at (-1, -1.8) {\small $x_2=3\,L/4$};
    \node [rotate=90] at (-1, -4.4) {\small $x_2=L/2$};
    \node [rotate=90] at (-1, -7) {\small $x_2=L/4$};
    \node [rotate=90] at (-1, -9.6) {\small $x_2=0$};
	\begin{groupplot}[group style={group size=5 by 5,horizontal sep=.5cm},ytick={0,5,11},width=.27\textwidth,xmin=-0.05, xmax=1.05, enlarge y limits = false, axis x line=bottom,axis y line=left,
    ]
	\def\myPlots{}
	\pgfplotsforeachungrouped \y in {4,3,2,1,0}
	{
		\pgfplotsforeachungrouped \x in {0,1,2,3,4}
		{
			\eappto\myPlots{%
				\noexpand\nextgroupplot
				\noexpand\addplot+[red,thick, mark=none]
				table[x index=1, y index=0, col sep=comma] {srcTikz/S33_sampling/Sig33_\x_\y.dat};
				\noexpand\addplot+[only marks,blue, mark=x,each nth point={9}]
				table[x index=2, y index=0, col sep=comma] {srcTikz/S33_sampling/Sig33_\x_\y.dat};			
			}
		}
	}
	\myPlots
	\end{groupplot}
	\end{tikzpicture}\tikzexternaldisable\else\includegraphics{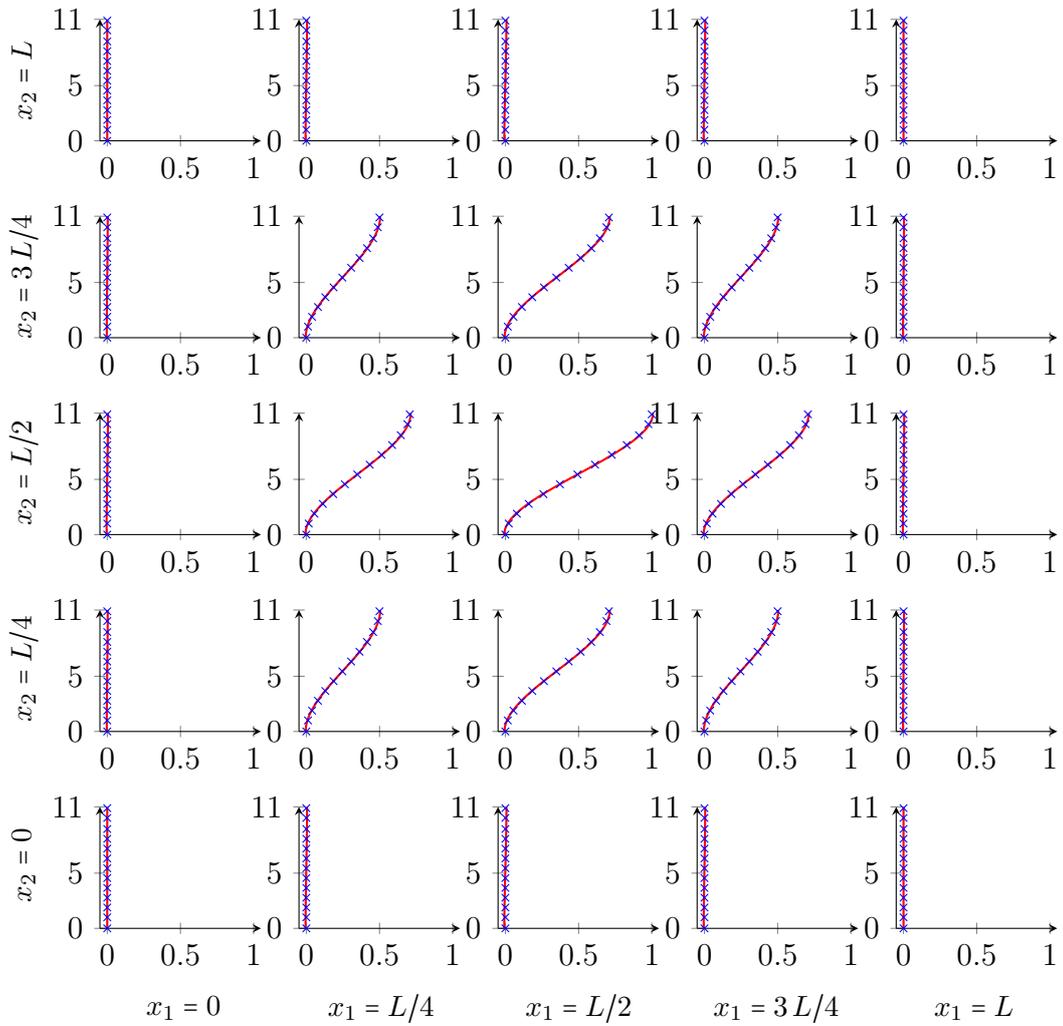}\fi
		\caption{Through-the-thickness $\bar{\sigma}_{33}$ profiles for several in plane sampling points.
		$L$ represents the total length of the plate, that for this case is  $L=220\,\text{mm}$ (being $L=S\,t$ with $t=11\,\text{mm}$ and $S=20$), while the number of layers is 11 (\textcolor{red}{$\boldsymbol{\leftrightline}$}post-processed solution, \textcolor{blue}{$\boldsymbol{\times}$} analytical solution~\cite{Pagano1970}).}
	\label{fig:samplingS33}
\end{figure}

\newpage
To show the effect of post-processing at different locations of the plate, in Figures \ref{fig:samplingS13}-\ref{fig:samplingS33} the out-of-plane stress state profile is recovered sampling the laminae every quarter of length in both in-plane directions, for the case of a length-to-thickness ratio equal to 20 and 11 layers.
\newpage
\subsection{Convergence behaviour}\label{subsec:conv}		
\begin{figure}[!htbp]
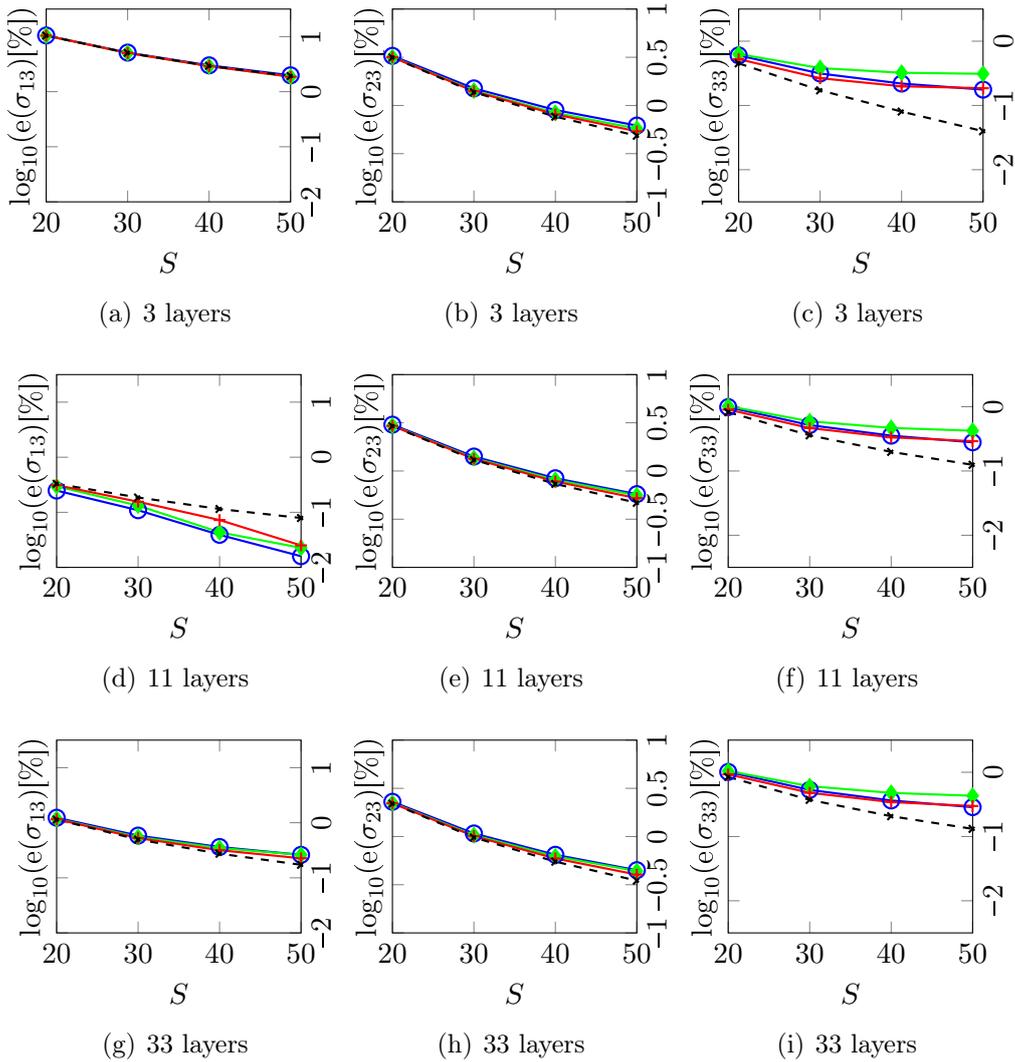

		\centering
	    \hspace{-7pt}
		\subfigure[3 layers \label{subfig-1:errors13_l3_664}]{\ifrecompiletikz\tikzsetnextfilename{fig_07_a}\tikzexternalenable\begin{tikzpicture}
\begin{axis}[width=.35\textwidth,
xlabel=$S$, 
ylabel={$\log_{10}$(e($\sigma_{13}$)[\%])},
xmin=20,xmax=50,
ymin=-2,ymax=1.5,
xtick={20,30,40,50},     
label style={font=\small},
tick label style={font=\small} ,
y label style={at={(axis description cs:0.39,.5)},anchor=south},
yticklabel pos=right,
yticklabel style={rotate=90},
]

\addplot+[color=blue,solid,mark=o,mark options={scale=1.5},thick] table[x index=0,y index=1, col sep=space] {srcTikz/error/664Error13_lay3el1.dat};
\addplot+[color=green,solid,mark=diamond*,mark options={scale=1.3},thick] table[x index=0,y index=1, col sep=space] {srcTikz/error/664Error13_lay3el2.dat};
\addplot+[color=red,solid,mark=+,mark options={scale=1.1},thick] table[x index=0,y index=1, col sep=space] {srcTikz/error/664Error13_lay3el4.dat};
\addplot+[color=black,solid,mark=x,mark options={scale=1},thick,dashed] table[x index=0,y index=1, col sep=space] {srcTikz/error/664Error13_lay3el8.dat};
\end{axis}
\end{tikzpicture}\tikzexternaldisable\else\includegraphics{fig_07_a.pdf}\fi}
		\hspace{-7pt}
		\subfigure[3 layers \label{subfig-2:errors23_l3_664}]{\ifrecompiletikz\tikzsetnextfilename{fig_07_b}\tikzexternalenable\begin{tikzpicture}
\begin{axis}[width=.35\textwidth,
xlabel=$S$, 
ylabel={$\log_{10}$(e($\sigma_{23}$)[\%])},
xmin=20,xmax=50,
ymin=-1,ymax=1,
xtick={20,30,40,50},
label style={font=\small},
tick label style={font=\small} ,
y label style={at={(axis description cs:0.39,.5)},anchor=south},
yticklabel pos=right,
yticklabel style={rotate=90},
]

\addplot+[color=blue,solid,mark=o,mark options={scale=1.5},thick] table[x index=0,y index=1, col sep=space] {srcTikz/error/664Error23_lay3el1.dat};
\addplot+[color=green,solid,mark=diamond*,mark options={scale=1.3},thick] table[x index=0,y index=1, col sep=space] {srcTikz/error/664Error23_lay3el2.dat};
\addplot+[color=red,solid,mark=+,mark options={scale=1.1},thick] table[x index=0,y index=1, col sep=space] {srcTikz/error/664Error23_lay3el4.dat};
\addplot+[color=black,solid,mark=x,mark options={scale=1},thick,dashed] table[x index=0,y index=1, col sep=space] {srcTikz/error/664Error23_lay3el8.dat};

\end{axis}
\end{tikzpicture}\tikzexternaldisable\else\includegraphics{fig_07_b.pdf}\fi}
		\hspace{-7pt}
		\subfigure[3 layers \label{subfig-3:errors33_l3_664}]{\ifrecompiletikz\tikzsetnextfilename{fig_07_c}\tikzexternalenable\begin{tikzpicture}
\begin{axis}[width=.35\textwidth,
    xmin=20,xmax=50,
ymin=-2.5,ymax=0.5,
    xtick={20,30,40,50},
xticklabels={$20$,$30$,$40$,$50$},
    y label style={at={(axis description cs:0.39,.5)},anchor=south},
    xlabel=$S$, 
    yticklabel pos=right,
    yticklabel style={rotate=90},
    ylabel={$\log_{10}$(e($\sigma_{33}$)[\%])},
    label style={font=\small},
    tick label style={font=\small} 
    ]

\addplot+[color=blue,solid,mark=o,mark options={scale=1.5},thick] table[x index=0,y index=1, col sep=space] {srcTikz/error/664Error33_lay3el1.dat};
\addplot+[color=green,solid,mark=diamond*,mark options={scale=1.3},thick] table[x index=0,y index=1, col sep=space] {srcTikz/error/664Error33_lay3el2.dat};
\addplot+[color=red,solid,mark=+,mark options={scale=1.1},thick] table[x index=0,y index=1, col sep=space] {srcTikz/error/664Error33_lay3el4.dat};
\addplot+[color=black,solid,mark=x,mark options={scale=1},thick,dashed] table[x index=0,y index=1, col sep=space] {srcTikz/error/664Error33_lay3el8.dat};

\end{axis}
\end{tikzpicture}\tikzexternaldisable\else\includegraphics{fig_07_c.pdf}\fi}\\
		\hspace{-7pt}
		\subfigure[11 layers \label{subfig-4:errors13_l11_664}]{\ifrecompiletikz\tikzsetnextfilename{fig_07_d}\tikzexternalenable\begin{tikzpicture}
\begin{axis}[width=.35\textwidth,
    xmin=20,xmax=50,
ymin=-2,ymax=1.5,
    xtick={20,30,40,50}, 
    y label style={at={(axis description cs:.39,.5)},anchor=south},
yticklabel pos=right,
    yticklabel style={rotate=90},
xlabel=$S$, 
ylabel={$\log_{10}$(e($\sigma_{13}$)[\%])},
    label style={font=\small},
tick label style={font=\small} 
   ]

\addplot+[color=blue,solid,mark=o,mark options={scale=1.5},thick] table[x index=0,y index=1, col sep=space] {srcTikz/error/664Error13_lay11el1.dat};
\addplot+[color=green,solid,mark=diamond*,mark options={scale=1.3},thick] table[x index=0,y index=1, col sep=space] {srcTikz/error/664Error13_lay11el2.dat};
\addplot+[color=red,solid,mark=+,mark options={scale=1.1},thick] table[x index=0,y index=1, col sep=space] {srcTikz/error/664Error13_lay11el4.dat};
\addplot+[color=black,solid,mark=x,mark options={scale=1},thick,dashed] table[x index=0,y index=1, col sep=space] {srcTikz/error/664Error13_lay11el8.dat};

\end{axis}
\end{tikzpicture}\tikzexternaldisable\else\includegraphics{fig_07_d.pdf}\fi}\hspace{-7pt}
		\subfigure[11 layers \label{subfig-5:errors23_l11_664}]{\ifrecompiletikz\tikzsetnextfilename{fig_07_e}\tikzexternalenable\begin{tikzpicture}
\begin{axis}[width=.35\textwidth,
xlabel=$S$, 
ylabel={$\log_{10}$(e($\sigma_{23}$)[\%])},
xmin=20,xmax=50,
ymin=-1,ymax=1,
xtick={20,30,40,50},
label style={font=\small},
tick label style={font=\small} ,
y label style={at={(axis description cs:0.39,.5)},anchor=south},
yticklabel pos=right,
yticklabel style={rotate=90},
]

\addplot+[color=blue,solid,mark=o,mark options={scale=1.5},thick] table[x index=0,y index=1, col sep=space] {srcTikz/error/664Error23_lay11el1.dat};
\addplot+[color=green,solid,mark=diamond*,mark options={scale=1.3},thick] table[x index=0,y index=1, col sep=space] {srcTikz/error/664Error23_lay11el2.dat};
\addplot+[color=red,solid,mark=+,mark options={scale=1.1},thick] table[x index=0,y index=1, col sep=space] {srcTikz/error/664Error23_lay11el4.dat};
\addplot+[color=black,solid,mark=x,mark options={scale=1},thick,dashed] table[x index=0,y index=1, col sep=space] {srcTikz/error/664Error23_lay11el8.dat};

\end{axis}
\end{tikzpicture}\tikzexternaldisable\else\includegraphics{fig_07_e.pdf}\fi}\hspace{-7pt}
		\subfigure[11 layers \label{subfig-6:errors33_l11_664}]{\ifrecompiletikz\tikzsetnextfilename{fig_07_f}\tikzexternalenable\begin{tikzpicture}
\begin{axis}[width=.35\textwidth,
xlabel=$S$, 
ylabel={$\log_{10}$(e($\sigma_{33}$)[\%])},
xmin=20,xmax=50,
ymin=-2.5,ymax=0.5,
xtick={20,30,40,50},
label style={font=\small},
tick label style={font=\small} ,
y label style={at={(axis description cs:0.39,.5)},anchor=south},
yticklabel pos=right,
yticklabel style={rotate=90},
]

\addplot+[color=blue,solid,mark=o,mark options={scale=1.5},thick] table[x index=0,y index=1, col sep=space] {srcTikz/error/664Error33_lay11el1.dat};
\addplot+[color=green,solid,mark=diamond*,mark options={scale=1.3},thick] table[x index=0,y index=1, col sep=space] {srcTikz/error/664Error33_lay11el2.dat};
\addplot+[color=red,solid,mark=+,mark options={scale=1.1},thick] table[x index=0,y index=1, col sep=space] {srcTikz/error/664Error33_lay11el4.dat};
\addplot+[color=black,solid,mark=x,mark options={scale=1},thick,dashed] table[x index=0,y index=1, col sep=space] {srcTikz/error/664Error33_lay11el8.dat};

\end{axis}
\end{tikzpicture}\tikzexternaldisable\else\includegraphics{fig_07_f.pdf}\fi}\\
		\hspace{-7pt}
		\subfigure[33 layers \label{subfig-7:errors13_l33_664}]{\ifrecompiletikz\tikzsetnextfilename{fig_07_g}\tikzexternalenable\begin{tikzpicture}
\begin{axis}[width=.35\textwidth,
xlabel=$S$, 
ylabel={$\log_{10}$(e($\sigma_{13}$)[\%])},
xmin=20,xmax=50,
ymin=-2,ymax=1.5,
xtick={20,30,40,50},       
label style={font=\small},
tick label style={font=\small} ,
y label style={at={(axis description cs:0.39,.5)},anchor=south},
yticklabel pos=right,
yticklabel style={rotate=90},
]

\addplot+[color=blue,solid,mark=o,mark options={scale=1.5},thick] table[x index=0,y index=1, col sep=space] {srcTikz/error/664Error13_lay33el1.dat};
\addplot+[color=green,solid,mark=diamond*,mark options={scale=1.3},thick] table[x index=0,y index=1, col sep=space] {srcTikz/error/664Error13_lay33el2.dat};
\addplot+[color=red,solid,mark=+,mark options={scale=1.1},thick] table[x index=0,y index=1, col sep=space] {srcTikz/error/664Error13_lay33el4.dat};
\addplot+[color=black,solid,mark=x,mark options={scale=1},thick,dashed] table[x index=0,y index=1, col sep=space] {srcTikz/error/664Error13_lay33el8.dat};

\end{axis}
\end{tikzpicture}\tikzexternaldisable\else\includegraphics{fig_07_g.pdf}\fi}\hspace{-7pt}
		\subfigure[33 layers \label{subfig-8:errors23_l33_664}]{\ifrecompiletikz\tikzsetnextfilename{fig_07_h}\tikzexternalenable\begin{tikzpicture}
\begin{axis}[width=.35\textwidth,
xlabel=$S$, 
ylabel={$\log_{10}$(e($\sigma_{23}$)[\%])},
xmin=20,xmax=50,
ymin=-1,ymax=1,
xtick={20,30,40,50},
label style={font=\small},
tick label style={font=\small} ,
y label style={at={(axis description cs:0.39,.5)},anchor=south},
yticklabel pos=right,
yticklabel style={rotate=90},
]

\addplot+[color=blue,solid,mark=o,mark options={scale=1.5},thick] table[x index=0,y index=1, col sep=space] {srcTikz/error/664Error23_lay33el1.dat};
\addplot+[color=green,solid,mark=diamond*,mark options={scale=1.3},thick] table[x index=0,y index=1, col sep=space] {srcTikz/error/664Error23_lay33el2.dat};
\addplot+[color=red,solid,mark=+,mark options={scale=1.1},thick] table[x index=0,y index=1, col sep=space] {srcTikz/error/664Error23_lay33el4.dat};
\addplot+[color=black,solid,mark=x,mark options={scale=1},thick,dashed] table[x index=0,y index=1, col sep=space] {srcTikz/error/664Error23_lay33el8.dat};

\end{axis}
\end{tikzpicture}\tikzexternaldisable\else\includegraphics{fig_07_h.pdf}\fi}\hspace{-7pt}
		\subfigure[33 layers \label{subfig-9:errors33_l33_664}]{\ifrecompiletikz\tikzsetnextfilename{fig_07_i}\tikzexternalenable\begin{tikzpicture}
\begin{axis}[width=.35\textwidth,
xlabel=$S$, 
ylabel={$\log_{10}$(e($\sigma_{33}$)[\%])},
xmin=20,xmax=50,
ymin=-2.5,ymax=0.5,
xtick={20,30,40,50},
label style={font=\small},
tick label style={font=\small} ,
y label style={at={(axis description cs:0.39,.5)},anchor=south},
yticklabel pos=right,
yticklabel style={rotate=90},
]

\addplot+[color=blue,solid,mark=o,mark options={scale=1.5},thick] table[x index=0,y index=1, col sep=space] {srcTikz/error/664Error33_lay33el1.dat};
\addplot+[color=green,solid,mark=diamond*,mark options={scale=1.3},thick] table[x index=0,y index=1, col sep=space] {srcTikz/error/664Error33_lay33el2.dat};
\addplot+[color=red,solid,mark=+,mark options={scale=1.1},thick] table[x index=0,y index=1, col sep=space] {srcTikz/error/664Error33_lay33el4.dat};
\addplot+[color=black,solid,mark=x,mark options={scale=1},thick,dashed] table[x index=0,y index=1, col sep=space] {srcTikz/error/664Error33_lay33el8.dat};

\end{axis}
\end{tikzpicture}\tikzexternaldisable\else\includegraphics{fig_07_i.pdf}\fi}
		\caption{Maximum relative percentage error evaluation at $x_1=x_2=0.25L$ for in-plane degree of approximation equal to 6 and out-of-plane degree of approximation equal to 4. Different length-to-thickness ratios $S$ are investigated for a number of layers equal to 3, 11, and 33 (\protect\markerone\hspace{0.2cm}1 knot span,\hspace{0.2cm}\protect\markertwo\hspace{0.2cm}2 knot spans,\hspace{0.2cm}\protect\markerthree\hspace{0.2cm}4 knot spans,\hspace{0.2cm}\protect\markerfour\hspace{0.2cm}8 knot spans).
		}
	    \label{fig:conv664}
\end{figure}

\begin{figure}[!htbp]
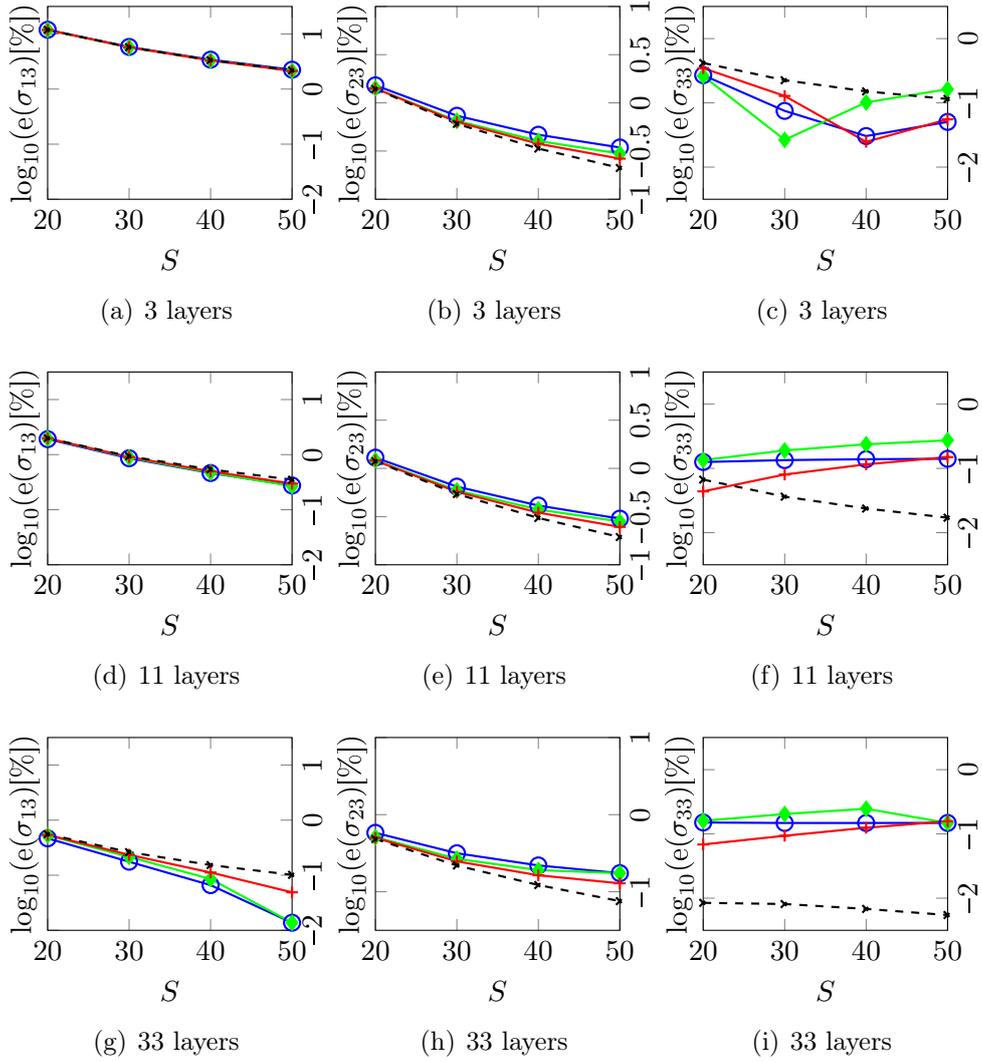

	\centering
	\hspace{-10pt}
	\subfigure[3 layers \label{subfig-1:errors13_l3_666}]{\ifrecompiletikz\tikzsetnextfilename{fig_08_a}\tikzexternalenable\begin{tikzpicture}
\begin{axis}[width=.35\textwidth,
xlabel=$S$, 
ylabel={$\log_{10}$(e($\sigma_{13}$)[\%])},
xmin=20,xmax=50,
ymin=-2,ymax=1.5,
xtick={20,30,40,50},
label style={font=\small},
tick label style={font=\small} ,
y label style={at={(axis description cs:0.39,.5)},anchor=south},
yticklabel pos=right,
yticklabel style={rotate=90},
]

\addplot+[color=blue,solid,mark=o,mark options={scale=1.5},thick] table[x index=0,y index=1, col sep=space] {srcTikz/error/666Error13_lay3el1.dat};
\addplot+[color=green,solid,mark=diamond*,mark options={scale=1.3},thick] table[x index=0,y index=1, col sep=space] {srcTikz/error/666Error13_lay3el2.dat};
\addplot+[color=red,solid,mark=+,mark options={scale=1.1},thick] table[x index=0,y index=1, col sep=space] {srcTikz/error/666Error13_lay3el4.dat};
\addplot+[color=black,solid,mark=x,mark options={scale=1},thick,dashed] table[x index=0,y index=1, col sep=space] {srcTikz/error/666Error13_lay3el8.dat};
\end{axis}
\end{tikzpicture}\tikzexternaldisable\else\includegraphics{fig_08_a.pdf}\fi}\hspace{-10pt}
	\subfigure[3 layers \label{subfig-4:errors23_l3_666}]{\ifrecompiletikz\tikzsetnextfilename{fig_08_b}\tikzexternalenable\begin{tikzpicture}
\begin{axis}[width=.35\textwidth,
xlabel=$S$, 
ylabel={$\log_{10}$(e($\sigma_{23}$)[\%])},
xmin=20,xmax=50,
ymin=-1,ymax=1,
xtick={20,30,40,50},
label style={font=\small},
tick label style={font=\small} ,
y label style={at={(axis description cs:0.39,.5)},anchor=south},
yticklabel pos=right,
yticklabel style={rotate=90},
]

\addplot+[color=blue,solid,mark=o,mark options={scale=1.5},thick] table[x index=0,y index=1, col sep=space] {srcTikz/error/666Error23_lay3el1.dat};
\addplot+[color=green,solid,mark=diamond*,mark options={scale=1.3},thick] table[x index=0,y index=1, col sep=space] {srcTikz/error/666Error23_lay3el2.dat};
\addplot+[color=red,solid,mark=+,mark options={scale=1.1},thick] table[x index=0,y index=1, col sep=space] {srcTikz/error/666Error23_lay3el4.dat};
\addplot+[color=black,solid,mark=x,mark options={scale=1},thick,dashed] table[x index=0,y index=1, col sep=space] {srcTikz/error/666Error23_lay3el8.dat};

\end{axis}
\end{tikzpicture}\tikzexternaldisable\else\includegraphics{fig_08_b.pdf}\fi}\hspace{-10pt}
	\subfigure[3 layers \label{subfig-7:errors33_l3_666}]{\ifrecompiletikz\tikzsetnextfilename{fig_08_c}\tikzexternalenable\begin{tikzpicture}
\begin{axis}[width=.35\textwidth,
    xmin=20,xmax=50,
ymin=-2.5,ymax=0.5,
    xtick={20,30,40,50},
xticklabels={$20$,$30$,$40$,$50$},
    y label style={at={(axis description cs:0.39,.5)},anchor=south},
    xlabel=$S$, 
    yticklabel pos=right,
    yticklabel style={rotate=90},
    ylabel={$\log_{10}$(e($\sigma_{33}$)[\%])},
    label style={font=\small},
    tick label style={font=\small} 
    ]

\addplot+[color=blue,solid,mark=o,mark options={scale=1.5},thick] table[x index=0,y index=1, col sep=space] {srcTikz/error/666Error33_lay3el1.dat};
\addplot+[color=green,solid,mark=diamond*,mark options={scale=1.3},thick] table[x index=0,y index=1, col sep=space] {srcTikz/error/666Error33_lay3el2.dat};
\addplot+[color=red,solid,mark=+,mark options={scale=1.1},thick] table[x index=0,y index=1, col sep=space] {srcTikz/error/666Error33_lay3el4.dat};
\addplot+[color=black,solid,mark=x,mark options={scale=1},thick,dashed] table[x index=0,y index=1, col sep=space] {srcTikz/error/666Error33_lay3el8.dat};

\end{axis}
\end{tikzpicture}\tikzexternaldisable\else\includegraphics{fig_08_c.pdf}\fi}\\
	\hspace{-10pt}
	\subfigure[11 layers \label{subfig-2:errors13_l11_666}]{\ifrecompiletikz\tikzsetnextfilename{fig_08_d}\tikzexternalenable\begin{tikzpicture}
\begin{axis}[width=.35\textwidth,
    xmin=20,xmax=50,
ymin=-2,ymax=1.5,
    xtick={20,30,40,50},
    y label style={at={(axis description cs:.39,.5)},anchor=south},
yticklabel pos=right,
    yticklabel style={rotate=90},
xlabel=$S$, 
ylabel={$\log_{10}$(e($\sigma_{13}$)[\%])},
    label style={font=\small},
tick label style={font=\small} 
   ]

\addplot+[color=blue,solid,mark=o,mark options={scale=1.5},thick] table[x index=0,y index=1, col sep=space] {srcTikz/error/666Error13_lay11el1.dat};
\addplot+[color=green,solid,mark=diamond*,mark options={scale=1.3},thick] table[x index=0,y index=1, col sep=space] {srcTikz/error/666Error13_lay11el2.dat};
\addplot+[color=red,solid,mark=+,mark options={scale=1.1},thick] table[x index=0,y index=1, col sep=space] {srcTikz/error/666Error13_lay11el4.dat};
\addplot+[color=black,solid,mark=x,mark options={scale=1},thick,dashed] table[x index=0,y index=1, col sep=space] {srcTikz/error/666Error13_lay11el8.dat};

\end{axis}
\end{tikzpicture}\tikzexternaldisable\else\includegraphics{fig_08_d.pdf}\fi}\hspace{-10pt}
	\subfigure[11 layers \label{subfig-5:errors23_l11_666}]{\ifrecompiletikz\tikzsetnextfilename{fig_08_e}\tikzexternalenable\begin{tikzpicture}
\begin{axis}[width=.35\textwidth,
xlabel=$S$, 
ylabel={$\log_{10}$(e($\sigma_{23}$)[\%])},
xmin=20,xmax=50,
ymin=-1,ymax=1,
xtick={20,30,40,50},
label style={font=\small},
tick label style={font=\small} ,
y label style={at={(axis description cs:0.39,.5)},anchor=south},
yticklabel pos=right,
yticklabel style={rotate=90},
]

\addplot+[color=blue,solid,mark=o,mark options={scale=1.5},thick] table[x index=0,y index=1, col sep=space] {srcTikz/error/666Error23_lay11el1.dat};
\addplot+[color=green,solid,mark=diamond*,mark options={scale=1.3},thick] table[x index=0,y index=1, col sep=space] {srcTikz/error/666Error23_lay11el2.dat};
\addplot+[color=red,solid,mark=+,mark options={scale=1.1},thick] table[x index=0,y index=1, col sep=space] {srcTikz/error/666Error23_lay11el4.dat};
\addplot+[color=black,solid,mark=x,mark options={scale=1},thick,dashed] table[x index=0,y index=1, col sep=space] {srcTikz/error/666Error23_lay11el8.dat};

\end{axis}
\end{tikzpicture}\tikzexternaldisable\else\includegraphics{fig_08_e.pdf}\fi}\hspace{-10pt}
	\subfigure[11 layers \label{subfig-8:errors33_l11_666}]{\ifrecompiletikz\tikzsetnextfilename{fig_08_f}\tikzexternalenable\begin{tikzpicture}
\begin{axis}[width=.35\textwidth,
xlabel=$S$, 
ylabel={$\log_{10}$(e($\sigma_{33}$)[\%])},
xmin=20,xmax=50,
ymin=-2.5,ymax=0.5,
xtick={20,30,40,50},
label style={font=\small},
tick label style={font=\small} ,
y label style={at={(axis description cs:0.39,.5)},anchor=south},
yticklabel pos=right,
yticklabel style={rotate=90},
]

\addplot+[color=blue,solid,mark=o,mark options={scale=1.5},thick] table[x index=0,y index=1, col sep=space] {srcTikz/error/666Error33_lay11el1.dat};
\addplot+[color=green,solid,mark=diamond*,mark options={scale=1.3},thick] table[x index=0,y index=1, col sep=space] {srcTikz/error/666Error33_lay11el2.dat};
\addplot+[color=red,solid,mark=+,mark options={scale=1.1},thick] table[x index=0,y index=1, col sep=space] {srcTikz/error/666Error33_lay11el4.dat};
\addplot+[color=black,solid,mark=x,mark options={scale=1},thick,dashed] table[x index=0,y index=1, col sep=space] {srcTikz/error/666Error33_lay11el8.dat};

\end{axis}
\end{tikzpicture}\tikzexternaldisable\else\includegraphics{fig_08_f.pdf}\fi}\\
	\hspace{-10pt}
	\subfigure[33 layers \label{subfig-3:errors13_l33_666}]{\ifrecompiletikz\tikzsetnextfilename{fig_08_g}\tikzexternalenable\begin{tikzpicture}
\begin{axis}[width=.35\textwidth,
xlabel=$S$, 
ylabel={$\log_{10}$(e($\sigma_{13}$)[\%])},
xmin=20,xmax=50,
ymin=-2,ymax=1.5,
xtick={20,30,40,50},
label style={font=\small},
tick label style={font=\small} ,
y label style={at={(axis description cs:0.39,.5)},anchor=south},
yticklabel pos=right,
yticklabel style={rotate=90},
]

\addplot+[color=blue,solid,mark=o,mark options={scale=1.5},thick] table[x index=0,y index=1, col sep=space] {srcTikz/error/666Error13_lay33el1.dat};
\addplot+[color=green,solid,mark=diamond*,mark options={scale=1.3},thick] table[x index=0,y index=1, col sep=space] {srcTikz/error/666Error13_lay33el2.dat};
\addplot+[color=red,solid,mark=+,mark options={scale=1.1},thick] table[x index=0,y index=1, col sep=space] {srcTikz/error/666Error13_lay33el4.dat};
\addplot+[color=black,solid,mark=x,mark options={scale=1},thick,dashed] table[x index=0,y index=1, col sep=space] {srcTikz/error/666Error13_lay33el8.dat};

\end{axis}
\end{tikzpicture}\tikzexternaldisable\else\includegraphics{fig_08_g.pdf}\fi}\hspace{-10pt}
	\subfigure[33 layers \label{subfig-6:errors23_l33_666}]{\ifrecompiletikz\tikzsetnextfilename{fig_08_h}\tikzexternalenable\begin{tikzpicture}
\begin{axis}[width=.35\textwidth,
xlabel=$S$, 
ylabel={$\log_{10}$(e($\sigma_{23}$)[\%])},
xmin=20,xmax=50,
ymin=-1.5,ymax=1,
xtick={20,30,40,50},
label style={font=\small},
tick label style={font=\small} ,
y label style={at={(axis description cs:0.39,.5)},anchor=south},
yticklabel pos=right,
yticklabel style={rotate=90},
]

\addplot+[color=blue,solid,mark=o,mark options={scale=1.5},thick] table[x index=0,y index=1, col sep=space] {srcTikz/error/666Error23_lay33el1.dat};
\addplot+[color=green,solid,mark=diamond*,mark options={scale=1.3},thick] table[x index=0,y index=1, col sep=space] {srcTikz/error/666Error23_lay33el2.dat};
\addplot+[color=red,solid,mark=+,mark options={scale=1.1},thick] table[x index=0,y index=1, col sep=space] {srcTikz/error/666Error23_lay33el4.dat};
\addplot+[color=black,solid,mark=x,mark options={scale=1},thick,dashed] table[x index=0,y index=1, col sep=space] {srcTikz/error/666Error23_lay33el8.dat};

\end{axis}
\end{tikzpicture}\tikzexternaldisable\else\includegraphics{fig_08_h.pdf}\fi}\hspace{-10pt}
	\subfigure[33 layers \label{subfig-9:errors33_l33_666}]{\ifrecompiletikz\tikzsetnextfilename{fig_08_i}\tikzexternalenable\begin{tikzpicture}
\begin{axis}[width=.35\textwidth,
xlabel=$S$, 
ylabel={$\log_{10}$(e($\sigma_{33}$)[\%])},
xmin=20,xmax=50,
ymin=-2.5,ymax=0.5,
xtick={20,30,40,50},
label style={font=\small},
tick label style={font=\small} ,
y label style={at={(axis description cs:0.39,.5)},anchor=south},
yticklabel pos=right,
yticklabel style={rotate=90},
]

\addplot+[color=blue,solid,mark=o,mark options={scale=1.5},thick] table[x index=0,y index=1, col sep=space] {srcTikz/error/666Error33_lay33el1.dat};
\addplot+[color=green,solid,mark=diamond*,mark options={scale=1.3},thick] table[x index=0,y index=1, col sep=space] {srcTikz/error/666Error33_lay33el2.dat};
\addplot+[color=red,solid,mark=+,mark options={scale=1.1},thick] table[x index=0,y index=1, col sep=space] {srcTikz/error/666Error33_lay33el4.dat};
\addplot+[color=black,solid,mark=x,mark options={scale=1},thick,dashed] table[x index=0,y index=1, col sep=space] {srcTikz/error/666Error33_lay33el8.dat};

\end{axis}
\end{tikzpicture}\tikzexternaldisable\else\includegraphics{fig_08_i.pdf}\fi}\\
	\caption{Maximum relative percentage error evaluation at $x_1=x_2=0.25L$ for degree of approximation equal to 6. Different length-to-thickness ratios $S$ are investigated for a number of layers equal to 3, 11, and 33 (\protect\markerone\hspace{0.2cm}1 knot span,\hspace{0.2cm}\protect\markertwo\hspace{0.2cm}2 knot spans,\hspace{0.2cm}\protect\markerthree\hspace{0.2cm}4 knot spans,\hspace{0.2cm}\protect\markerfour\hspace{0.2cm}8 knot spans).}
	\label{fig:conv666}
\end{figure}

 In order to validate the proposed approach in a wider variety of cases, computations with a
 different ratio between the thickness of the plate and its length are performed respectively
 for 3, 11, and 33 layers, considering an increasing number of knot spans. Figures~\ref{fig:conv664} and~\ref{fig:conv666} assess the convergence behaviour of the method, adopting the following error definition
	\begin{equation}
		\text{e}(\sigma_{ij})=\frac{\max(|\sigma_{ij}^\text{analytic}-\sigma_{ij}^\text{recovered}|)}{\max(|\sigma_{ij}^\text{analytic}|)}\,.\label{eq:error}
	\end{equation}
Note that relation \eqref{eq:error} is used only to estimate the error inside the domain to avoid indeterminate forms.
Different combinations of degree of approximations have been also considered. A poorer out-of-plane stress approximation is obtained using a degree equal to 4 in every direction, and, in addition, with this choice locking phenomena may occur for increasing values of length-to-thickness ratio. Therefore, we conclude that using a degree of approximation equal to 6  in-plane and equal to 4 through the thickness seems to be a reasonable choice to correctly reproduce the 3D stress state. Using instead uniform approximation degrees $p=q=r=6$ does not seem to significantly improve the results (see Figures~\ref{fig:conv664},~\ref{fig:conv666}, and Table~\ref{tab:errorl11}).
The post-processing method provides better results for increasing values of length-to-thickness ratio and number of layers and therefore proves to be particularly convenient for very large and thin plates. 
This is clear since a laminae with a large number of thin layers resembles a plate with average properties. 
What really stands out looking at the displayed mesh sensitivity results, is the fact that collocation perfectly captures the plates behaviour not only using one element through the thickness but also employing only one knot span in the plane of the plate. A single element of degrees $p=q=6$ and $r=4$, comprising a total of 7x7x5 collocation points, is able to provide for this example maximum percentage errors of 5\% or lower (and of 1\% or lower in the cases of 11 and 33 layers) for $S=30$ or larger.

\newpage
Quantitative results are presented in Table~\ref{tab:errorl11} for various plate cases, considering a number of layers equal to 11 and 10 collocation points for each in-plane parametric direction. Different number of layers (i.e., 3 and 33) are instead investigated in Appendix~\ref{sec:appendix}. Increasing length-to-thickness ratios, namely 20, 30, 40, and 50 are considered and the maximum relative error results is reported for a reference point located at $x_1=x_2=0.25L$.
Also different degrees of approximation are investigated. Given these results, we conclude that using an out-of-plane degree of approximation equal to 4 leads to a sufficiently accurate stress state. Furthermore the out-of-plane stress profile reconstruction shows a remarkable improvement for increasing values of number of layers and slenderness parameter $S$.
\begin{table}[!htbp]\caption{Simply supported composite plate under sinusoidal load with a number of layers equal to 11. Out-of-plane stress state maximum relative error with respect to Pagano's solution~\cite{Pagano1970} at $x_1=x_2=0.25L$. Comparing the isogeometric collocation-based homogenized single element approach (IGA-C) and the coupled post-processing technique (IGA-C+PP) for different approximation degrees.} 	\vspace{0.5cm}
		\centering	
		\begin{tabular}{? l | c ? c | c | c ? c | c | c ?}
			\thickhline
			\multicolumn{2}{? c ?}{\textbf{Degree}}
			& \multicolumn{3}{c ?}
	{$p=q=6$, $r=4$}&\multicolumn{3}{c ?}{$p=q=r=6$}\Tstrut\Bstrut\\\thickhline
\multirow{2}{*}{\textbf{S}} &\multirow{2}{*}{\textbf{Method}}
&$e(\sigma_{13})$ 
&$e(\sigma_{23})$ 
&$e(\sigma_{33})$ 
&$e(\sigma_{13})$ 
&$e(\sigma_{23})$ 
&$e(\sigma_{33})$  
			\Tstrut\Bstrut\\\cline{3-8}
			& &[\%]&[\%]&[\%]&[\%]&[\%]&[\%]
			\Tstrut\Bstrut\\\thickhline
			\multirow{2}{*}{20} & IGA-C & 97.6 & 56.7 & 6.34 & 96.6 & 56.1 & 6.31\Tstrut\Bstrut\\
			& IGA-C+PP & 0.31 & 2.94 & 0.90 & 1.97 & 1.20 & 0.05\Tstrut\Bstrut\\\hline
			\multirow{2}{*}{30} & IGA-C & 98.7 & 55.6 & 6.36 & 98.3 & 55.4 & 6.34\Tstrut\Bstrut\\
			& IGA-C+PP & 0.16 & 1.34 & 0.47 & 0.91 & 0.57 & 0.08\Tstrut\Bstrut\\\hline
			\multirow{2}{*}{40} & IGA-C & 99.2 & 55.3 & 6.37 & 98.9 & 55.2 & 6.36\Tstrut\Bstrut\\
			& IGA-C+PP & 0.07 & 0.78 & 0.34 & 0.50 & 0.35 & 0.12\Tstrut\Bstrut\\\hline
			\multirow{2}{*}{50} & IGA-C & 99.4 & 55.1 & 6.38 & 99.2 & 55.1 & 6.38\Tstrut\Bstrut\\
			& IGA-C+PP & 0.03 & 0.52 & 0.29 & 0.30 & 0.25 & 0.15\Tstrut\Bstrut\\\thickhline
		\end{tabular}\label{tab:errorl11}
\end{table}
\newpage
\section{Conclusions}\label{sec:conclusions}
In this paper we present a new approach to simulate laminated plates characterized by a symmetric distribution of plies. This technique combines a 3D collocation isogeometric analysis with a post-processing step procedure based on equilibrium equations.
Since we adopt a single-element appoach, to take into account variation
through the plate thickness of the material properties, we average the constitutive behaviour of each layer considering an homogeneized response.
Following this simple approach, we showed that acceptable results can be obtained only in terms of displacements and in-plain stresses. Therefore, we propose to perform a post-processing step which requires the shape functions to be highly continuous.
This continuity demand is fully granted by typical IGA shape functions. After the post-processing correction is applied, good results are recovered also in terms of out-of-plane stresses, even for very coarse meshes.
The post-processing stress-recovery technique is only based on the integration through the thickness of equilibrium equations,
and all the required components can be easily computed differentiating the displacement solution. 
Several numerical tests are carried out to test the sensitivity of the proposed technique to different length-to-thickness ratios and number
of layers. Regardless of the number of layers, the
method gives better results the thinner the composites are.
Multiple numbers of alternated layers and sequence of stacks (both even and odd) have been studied in our applications. Neverthless only tests which consider an odd number of layers or an odd disposition of an even number of stacks show good results as expected because considering a homogenized response of the material is effective only for symmetric distributions of plies.
Further research topics currently under investigation consist in the extension of this approach to
more complex problems involving curved geometries and large deformations.
\section*{Acknowledgments}\label{sec:acknowl}
This work was partially supported by Fondazione Cariplo -- Regione Lombardia through the project "Verso nuovi strumenti di simulazione super veloci ed accurati basati sull'analisi isogeometrica", within the program RST -- rafforzamento.\\
\indent P. Antolin was partially supported by the European Research council through the H2020 ERC Advanced Grant 2015 n.694515 CHANGE.
\begin{appendices}
\section{}
\label{sec:appendix}
Results in terms of maximum relative error considering a plate with a number of layers equal to 3 and 33 are herein presented for a reference point located at $x_1=x_2=0.25L$.
Increasing length-to-thickness ratios, namely 20, 30, 40, and 50 are investigated for different degrees of approximations (i.e., $p=q=6$ and $r=4$, and $p=q=r=6$), using 10 collocation points for each in-plane parametric direction and one element through-the-thickness.

\begin{table}[!htbp]\caption{Simply supported composite plate under sinusoidal load with a number of layers equal to 3. Out-of-plane stress state maximum relative error with respect to Pagano's solution~\cite{Pagano1970} at $x_1=x_2=0.25L$. Comparing the isogeometric collocation-based homogenized single element approach (IGA-C) and the coupled post-processing technique (IGA-C+PP) for different approximation degrees.} 	\vspace{0.5cm}
	\centering	
		\begin{tabular}{? l | c ? c | c | c ? c | c | c ?}
	\thickhline
	\multicolumn{2}{? c ?}{\textbf{Degree}}
	& \multicolumn{3}{c ?}
	{$p=q=6$, $r=4$}&\multicolumn{3}{c ?}{$p=q=r=6$}\Tstrut\Bstrut\\\thickhline
	\multirow{2}{*}{\textbf{S}} &\multirow{2}{*}{\textbf{Method}}
	&$e(\sigma_{13})$ 
	&$e(\sigma_{23})$ 
	&$e(\sigma_{33})$ 
	&$e(\sigma_{13})$ 
	&$e(\sigma_{23})$ 
	&$e(\sigma_{33})$ 
	\Tstrut\Bstrut\\\cline{3-8}
	& &[\%]&[\%]&[\%]&[\%]&[\%]&[\%]
	\Tstrut\Bstrut\\\thickhline
	\multirow{2}{*}{20} & IGA-C & 292 & 57.2 & 5.80& 291 & 57.2 & 5.79\Tstrut\Bstrut\\
	& IGA-C+PP & 10.4 & 3.16 & 0.54 & 11.9 & 1.41 & 0.33\Tstrut\Bstrut\\\hline 
	\multirow{2}{*}{30} & IGA-C & 311 & 57.5 & 5.77& 311 & 57.5 & 5.77\Tstrut\Bstrut\\
	& IGA-C+PP & 5.05 & 1.40 & 0.28 & 5.75 & 0.63 & 0.11\Tstrut\Bstrut\\\hline
	\multirow{2}{*}{40} & IGA-C & 319 & 57.6 & 5.76 & 319 & 57.6 & 5.76\Tstrut\Bstrut\\
	& IGA-C+PP & 2.91 & 0.81 & 0.21 & 3.32 & 0.38 & 0.02\Tstrut\Bstrut\\\hline
	\multirow{2}{*}{50} & IGA-C & 323 & 57.6 & 5.76& 322 & 57.6 & 5.76\Tstrut\Bstrut\\
	& IGA-C+PP & 1.87 & 0.54 & 0.20 & 2.14 & 0.26 & 0.07\Tstrut\Bstrut\\\thickhline
	\end{tabular}\label{tab:errorl3}
\end{table}

\begin{table}[!htbp]\caption{Simply supported composite plate under sinusoidal load with a number of layers equal to 33. Out-of-plane stress state maximum relative error with respect to Pagano's solution~\cite{Pagano1970} at $x_1=x_2=0.25L$. Comparing the isogeometric collocation-based homogenized single element approach (IGA-C) and the coupled post-processing technique (IGA-C+PP) for different approximation degrees.} 	\vspace{0.5cm}
	\centering	
		\begin{tabular}{? l | c ? c | c | c ? c | c | c ?}
	\thickhline
	\multicolumn{2}{? c ?}{\textbf{Degree}}
	& \multicolumn{3}{c ?}
	{$p=q=6$, $r=4$}&\multicolumn{3}{c ?}{$p=q=r=6$}\Tstrut\Bstrut\\\thickhline
\multirow{2}{*}{\textbf{S}} &\multirow{2}{*}{\textbf{Method}}
&$e(\sigma_{13})$ 
&$e(\sigma_{23})$ 
&$e(\sigma_{33})$ 
&$e(\sigma_{13})$ 
&$e(\sigma_{23})$ 
&$e(\sigma_{33})$  
	\Tstrut\Bstrut\\\cline{3-8}
	& &[\%]&[\%]&[\%]&[\%]&[\%]&[\%]
	\Tstrut\Bstrut\\\thickhline
	\multirow{2}{*}{20} & IGA-C & 81.6 & 69.7 & 6.33 & 80.7 & 68.9 & 6.33\Tstrut\Bstrut\\
	& IGA-C+PP & 1.16 & 2.21 & 0.93 & 0.54 & 0.50 & 0.07\Tstrut\Bstrut\\\hline
	\multirow{2}{*}{30} & IGA-C & 81.5 & 69.0 & 6.34 & 81.2 & 68.7 & 6.34\Tstrut\Bstrut\\
	& IGA-C+PP & 0.53 & 1.01 & 0.48 & 0.23 & 0.25 & 0.09\Tstrut\Bstrut\\\hline
	\multirow{2}{*}{40} & IGA-C & 81.5 & 68.7 & 6.35& 81.3 & 68.6 & 6.34\Tstrut\Bstrut\\
	& IGA-C+PP & 0.32 & 0.59 & 0.34 & 0.11 & 0.16 & 0.12\Tstrut\Bstrut\\\hline
	\multirow{2}{*}{50} & IGA-C & 81.6 & 68.6 & 6.35 & 81.4 & 68.5 & 6.35\Tstrut\Bstrut\\
	& IGA-C+PP & 0.23 & 0.40 & 0.30 & 0.05 & 0.13 & 0.16\Tstrut\Bstrut\\\thickhline
	\end{tabular}\label{tab:errorl33}
\end{table}
\end{appendices}

\newpage
\section*{}
\bibliographystyle{unsrt}   

\end{document}